\input amstex
\documentstyle{amsppt}
\NoBlackBoxes
%
\chardef\tempcat=\the\catcode`\@
\catcode`\@=11

\font\fourteenbf=cmbx10 scaled \magstep2

\def\QED{\nobreak\quad\ifmmode\roman{Q.E.D.}\else{\rm Q.E.D.}\fi}
\def\keyedby#1{}
\def\affil#1\endaffil{}
\def\ratitle{\setbox0=\hbox{\fourteenbf R}
                 \advance\baselineskip by \ht0  
                 \centerline{\fourteenbf RESEARCH ANNOUNCEMENTS}%
                 \baselineskip45pt}     
\def\date#1\enddate{\thanks Received by the editors #1 \endthanks}
\def\shorttitle#1{\rightheadtext{#1}}
\def\shortauthor#1{\leftheadtext{#1}}
\def\subjclassrev #1{\thanks {{\rm1980 {\it Mathematics Subject
       Classification} (1985 {\it Revision}). #1}}\endthanks}
\def\ml #1\endml{\email #1\endemail}
\long\def\ext #1\endext{\block #1\endblock}
\def\ac #1\endac{{\monograph@true \indenti=0pc 
                 \heading Contents\endheading \toc #1\endtoc}}
\def\ch #1\endch{\def\subheadfont@{\it} \subheading{#1} \def\subheadfont@{\bf}}
\let\thm=\proclaim
\let\ethm=\endproclaim

\let\rem=\remark
\let\endrem=\endremark

\def\RM#1{{\rm#1}}
\def\<{\nobreak\kern-\mathsurround\nobreak}
\everydisplay={\def\.{\thinspace.}}
\def\fighere{\def\figure@type{3}\ins@}

\loadbold 

\def\cprime{\/{\m@th$'$}}

{\let\enddocument=\relax
\gdef\bookrev{%
	\long\gdef\revtop ##1\endrevtop{\def\revtop@{{\interlinepenalty10000
		\tenpoint
		\everypar{\hangindent=\parindent}
                \noindent ##1\endgraf}}}%

	\gdef\reviewer ##1\endreviewer{\def\revr@{{\leftskip0pt plus1fil
	  \rightskip0pt \parfillskip0pt \tenpoint\smc ##1\endgraf}}}%

	\gdef\affil ##1\endaffil{\def\affil@{{\leftskip0pt plus1fil
	  \rightskip0pt \parfillskip0pt \tenpoint\smc ##1\endgraf}}}%

      \gdef\bookrevheading{\hrule height\z@
                 \vskip-\topskip
                 \setbox0=\hbox{\fourteenbf B}
                 \advance\baselineskip by \ht0  
                 \centerline{\fourteenbf BOOK REVIEW}%
                 \baselineskip45pt\hbox{}\baselineskip0pt}

   \gdef\document{%
        \bookrevheading
        \gdef\bookrevheading{\penalty-100\vskip2pc plus2pc}
   	\leftheadtext{BOOK REVIEW}%
	\rightheadtext{BOOK REVIEW}%
	\nobreak
	\vskip8pt  plus5pt minus2pt
	\nobreak
	\revtop@
	\nobreak
	\vskip12pt plus5pt minus2pt
	\tenpoint }%

	\outer\def\enddocument{\nobreak
	\vskip6pt minus6pt
	\nobreak
	\revr@
        \nobreak
        \affil@
	\end }%

}
}

\def\prelogo{}


\def\jourlogo{\hbox{\vbox to0pt{\prelogo
   \vbox{\sixrm \baselineskip 6.5pt
   \parindent 0pt APPEARED IN BULLETIN OF THE\hfil\break
   AMERICAN MATHEMATICAL SOCIETY\hfil\break
   Volume \cvol@, Number \cvolno@, \cmonth@\ \cvolyear@, Pages \cpgs@\endgraf}%
   \vss}}}


\def\cvol#1{\gdef\cvol@{\ignorespaces#1\unskip}}
\def\cvolno#1{\gdef\cvolno@{\ignorespaces#1\unskip}}
\def\cmonth#1{\gdef\cmonth@{\ignorespaces#1\unskip}}
\def\cvolyear#1{\gdef\cvolyear@{\ignorespaces#1\unskip}}
\def\cyear#1{\gdef\cyear@{\ignorespaces#1\unskip}\cyear@@#100000\end@}
\def\cpgs#1{\gdef\cpgs@{\ignorespaces#1\unskip}}

\def\cyear@@#1#2#3#4#5\end@{\gdef\cyearmodc@{#3#4}%
        \gdef\cyearmodcHold@{#3#4}}

\cvol{000}
\cvolno{0}
\cmonth{}
\cyear{0000}
\cvolyear{0000}
\cpgs{}

\font\sixsy=cmsy6

\def\copyrightline@{\baselineskip1.75pc
    \rightline{%
        \vbox{\sixrm \textfont2=\sixsy \baselineskip 7pt
            \halign{\hfil##\cr
                \copyright\cyear@\ American Mathematical Society\cr
                 0273-0979/\cyearmodc@\ \$1.00 + \$.25 per page\cr }}}}

\def\cyearmodc#1{\gdef\cyearmodc@{\ignorespaces#1\unskip}}


\let\logo@=\copyrightline@
\def\pretitle{\jourlogo \vskip24\p@ plus12\p@ minus12\p@}



\def\shoveright#1{\omit\span\omit\span\omit 
  \hfil$\@lign\displaystyle{{}#1}\m@th$%
  \iftagsleft@ 
    \ifx\undefined\displaywidth@ 
    \else\kern-\displaywidth@\fi
  \fi}


\def\aligned@{\bgroup\vspace@\Let@
 \def\shoveright##1{\omit\span\egroup\hfill\bgroup##1}%
 \ifinany@\else\openup\jot\fi\ialign
 \bgroup\hfil\strut@$\m@th\displaystyle{##}$&
 $\m@th\displaystyle{{}##}$\hfil\crcr}


\def\gathered{\null\,\vcenter\bgroup\vspace@\Let@
 \def\shoveright##1{##1\aftergroup\aftergroup\aftergroup\hfilneg}%
 \ifinany@\else\openup\jot\fi\ialign
 \bgroup\hfil\strut@$\m@th\displaystyle{##}$\hfil\crcr}


\def\gather{\RIfMIfI@\DN@{\onlydmatherr@\gather}\else
 \ingather@true\inany@true\def\tag{&}%
 \vspace@\allowdisplaybreak@\displaybreak@\intertext@
 \displ@y\Let@
 \def\shoveright##1{\omit\span\omit
   \hfil\llap{\strut@$\m@th\displaystyle{##1}$}%
   \iftagsleft@
     \ifx\undefined\gdisplaywidth@ 
     \else\kern-\gdisplaywidth@\fi
   \fi}%
 \iftagsleft@\DN@{\csname gather \endcsname}\else
  \DN@{\csname gather \space\endcsname}\fi\fi
 \else\DN@{\onlydmatherr@\gather}\fi\next@}

%
%
\newskip\colspacing      \colspacing=1em
%
%
\def\varquad{\hskip\colspacing\relax}
\def\reduce#1{\advance#1 by -}
%
%
\def\halfspace{\hskip .5\colspacing\relax}
%
%
\def\ruleh{\noalign{\hrule}}
%
%
%

%
\def\t@bleside{\ifvmode\hrule \else\vrule\fi\relax}
\newif\iftop@ 
\def\toptablecaption{\top@true\tablecaption}
\def\tablecapf@nt{\tenpoint}
%
%
\def\tablecaption#1{\iftop@\else\medskip\fi\begingroup\tablecapf@nt\Let@
  \centerline{\vbox{\ialign{\hfil##\hfil\cr#1\crcr}}}%
  \endgroup\iftop@\medskip\fi}
%
%
\def\colheading#1{\vbox{\normalbaselines\Let@\ialign{\rm\hfil##\hfil\cr
        \tablestrut depth\z@ #1\unskip\tablestrut height\z@\crcr}}}
%
%
\def\table{\relaxnext@\DN@{\ifx\next\nobox\def\t@bleside{}\fi
 \relaxnext@{}\,\t@bleside\vcenter\bgroup\Let@\vspace@
 \offinterlineskip\t@bleside\ialign\bgroup
  \halfspace\hfil$\tablestrut####$\hfil&&\varquad\hfil$####$\hfil\crcr
}\futurelet\next\next@}
%
\def\nobox{}
\def\endtable{\crcr\egroup\t@bleside\egroup\t@bleside\,}
%

\def\matrix{\null\,\vcenter\bgroup\Let@\vspace@
 \def\ruleh{\noalign{\vskip\tw@\p@\hrule\vskip\thr@@\p@}}%
 \def\everyr@w{}%
 \normalbaselines\openup\spreadmlines@\m@th\ialign
 \bgroup\hfil$##$\hfil&&\varquad\hfil$##$\hfil\crcr
 \Mathstrut@\crcr\noalign{\kern-\baselineskip}}

\def\everyr@w{\halfspace\tablestrut}

\newdimen\rowspacing  


\def\normalbaselines{\lineskip\normallineskip
  \baselineskip\normalbaselineskip \lineskiplimit\normallineskiplimit
  \rowspacing1.6\normalbaselineskip}

%
%
\def\tablestrut{\vrule width\z@ height.67\rowspacing
  depth.33\rowspacing}

\def\deeper#1{\begingroup
  \dimen@.33\rowspacing\advance\dimen@#1\tablestrut depth\dimen@
  \relax\endgroup}
\def\higher#1{\begingroup
  \dimen@.67\rowspacing\advance\dimen@ #1\tablestrut height\dimen@
  \relax\endgroup}

\def\format@#1\\{\def\preamble@{#1}%
 \def\c{\hfil$\the\hashtoks@$\hfil}%
 \def\r{\hfil$\the\hashtoks@$}%
 \def\l{$\the\hashtoks@$\hfil}%
 \def\ctext{\hfil\the\hashtoks@\hfil}%
 \def\rtext{\hfil\the\hashtoks@}%
 \def\ltext{\the\hashtoks@\hfil}%
 \setboxz@h{\xdef\Preamble@{\everyr@w\preamble@}}%
 \ifnum`{=0 \fi\iffalse}\fi
 \ialign\bgroup\span\Preamble@\crcr}

\catcode`\@=\tempcat


\keyedby{BULL225E/PAZ}

\define\a{\alpha}

\define\n{\noindent}

\define\EPS{\varepsilon}
\define\e{\enspace}
\define\th{\thinspace}

\define\di{\displaystyle }
\define\boxit#1{\vbox{\hrule\hbox{\vrule\kern1pt\vbox{%
\kern1pt#1\kern1pt}
\kern1pt\vrule}\hrule}}

\define\va{\varphi}

\define\1{\Omega}

\define\2{\Delta}
\define\RA{\rightarrow}


\define\de{\delta}

\define\rn{\Bbb R^n}




\define\il{\int}
\define\ur{u_\rho }


\define\SS{\Bigl|\Bigr|}
\define\3{\bigl|\bigr|}
\define\dt{{\delta}}
\define\lra{\longrightarrow}

\define\bbbr{\Bbb{R}}  
\define\bbbn{\Bbb{N}}  

\define\bbbc{{\mathchoice 
{\setbox0=\hbox{$\displaystyle\rm C$}\hbox
{\hbox
to0pt{\kern0.4\wd0\vrule height0.9\ht0\hss}\box0}}
{\setbox0=\hbox{$\textstyle\rm C$}\hbox{\hbox
to0pt{\kern0.4\wd0\vrule height0.9\ht0\hss}\box0}}
{\setbox0=\hbox{$\scriptstyle\rm C$}\hbox{\hbox
to0pt{\kern0.4\wd0\vrule height0.9\ht0\hss}\box0}}
{\setbox0=\hbox{$\scriptscriptstyle\rm C$}\hbox{\hbox
to0pt{\kern0.4\wd0\vrule height0.9\ht0\hss}\box0}}}}
\define\bbbe{{\mathchoice {\setbox0=\hbox{\smalletextfont 
e}
\hbox{\raise
0.1\ht0\hbox to0pt{\kern0.4\wd0\vrule width0.3pt
height0.7\ht0\hss}\box0}}
{\setbox0=\hbox{\smalletextfont e}\hbox{\raise 0.1\ht0\hbox
to0pt{\kern0.4\wd0\vrule width0.3pt 
height0.7\ht0\hss}\box0}}
{\setbox0=\hbox{\smallescriptfont e}\hbox{\raise 
0.1\ht0\hbox
to0pt{\kern0.5\wd0\vrule width0.2pt 
height0.7\ht0\hss}\box0}}
{\setbox0=\hbox{\smallescriptscriptfont e}\hbox{\raise
0.1\ht0\hbox to0pt{\kern0.4\wd0\vrule width0.2pt
height0.7\ht0\hss}\box0}}}}
\define\bbbq{{\mathchoice 
{\setbox0=\hbox{$\displaystyle\rm Q$}
\hbox{\raise
0.15\ht0\hbox to0pt{\kern0.4\wd0\vrule 
height0.8\ht0\hss}\box0}}
{\setbox0=\hbox{$\textstyle\rm Q$}\hbox{\raise
0.15\ht0\hbox to0pt{\kern0.4\wd0\vrule 
height0.8\ht0\hss}\box0}}
{\setbox0=\hbox{$\scriptstyle\rm Q$}\hbox{\raise
0.15\ht0\hbox to0pt{\kern0.4\wd0\vrule 
height0.7\ht0\hss}\box0}}
{\setbox0=\hbox{$\scriptscriptstyle\rm Q$}\hbox{\raise
0.15\ht0\hbox to0pt{\kern0.4\wd0\vrule 
height0.7\ht0\hss}\box0}}}}
\define\bbbt{{\mathchoice {\setbox0=\hbox{$\displaystyle\rm
T$}\hbox{\hbox to0pt{\kern0.3\wd0\vrule 
height0.9\ht0\hss}\box0}}
{\setbox0=\hbox{$\textstyle\rm T$}\hbox{\hbox
to0pt{\kern0.3\wd0\vrule height0.9\ht0\hss}\box0}}
{\setbox0=\hbox{$\scriptstyle\rm T$}\hbox{\hbox
to0pt{\kern0.3\wd0\vrule height0.9\ht0\hss}\box0}}
{\setbox0=\hbox{$\scriptscriptstyle\rm T$}\hbox{\hbox
to0pt{\kern0.3\wd0\vrule height0.9\ht0\hss}\box0}}}}
\define\bbbs{{\mathchoice
{\setbox0=\hbox{$\displaystyle     \rm 
S$}\hbox{\raise0.5\ht0\hbox
to0pt{\kern0.35\wd0\vrule height0.45\ht0\hss}\hbox
to0pt{\kern0.55\wd0\vrule height0.5\ht0\hss}\box0}}
{\setbox0=\hbox{$\textstyle        \rm 
S$}\hbox{\raise0.5\ht0\hbox
to0pt{\kern0.35\wd0\vrule height0.45\ht0\hss}\hbox
to0pt{\kern0.55\wd0\vrule height0.5\ht0\hss}\box0}}
{\setbox0=\hbox{$\scriptstyle      \rm 
S$}\hbox{\raise0.5\ht0\hbox
to0pt{\kern0.35\wd0\vrule 
height0.45\ht0\hss}\raise0.05\ht0\hbox
to0pt{\kern0.5\wd0\vrule height0.45\ht0\hss}\box0}}
{\setbox0=\hbox{$\scriptscriptstyle\rm 
S$}\hbox{\raise0.5\ht0\hbox
to0pt{\kern0.4\wd0\vrule 
height0.45\ht0\hss}\raise0.05\ht0\hbox
to0pt{\kern0.55\wd0\vrule height0.45\ht0\hss}\box0}}}}
\define\bbbz{{\mathchoice {\hbox{$\sans\textstyle 
Z\kern-0.4em Z$}}
{\hbox{$\sans\textstyle Z\kern-0.4em Z$}}
{\hbox{$\sans\scriptstyle Z\kern-0.3em Z$}}
{\hbox{$\sans\scriptscriptstyle Z\kern-0.2em Z$}}}}
\topmatter
\cvol{26}
\cvolyear{1992}
\cmonth{Jan}
\cyear{1992}
\cvolno{1}
\cpgs{53-86}
\title Semilinear wave equations\endtitle
\shorttitle{Semilinear wave equations}
\author Michael Struwe\endauthor
\shortauthor{Michael Struwe}
\address Mathematik, ETH-Zentrum, CH-8092 Z\"urich,
Switzerland\endaddress
\date July 11, 1990 and, in revised form, December 7, 
1990\enddate
\subjclassrev{Primary 35L05, 35A05, 35-02}
\abstract We survey existence and regularity results for 
semi-linear wave equations. In
particular, we review the recent regularity results for the
$u^5$-Klein Gordon equation by Grillakis and this author 
and give a
self-contained, slightly simplified proof.\endabstract
\endtopmatter

\document
\ac
\heading 1. Introduction\endheading  
\heading 2. Preliminaries\endheading 
\heading 3. Rauch's result\endheading 
\heading 4. Large data\endheading 
\heading 5. A remark on the super-critical case\endheading
\endac

\heading 1. Introduction\endheading 
In this survey we shall be interested in  initial value 
problems for nonlinear wave
equations of the type 
$$u_{tt}-\2 u+g(u)=0\quad\text{{\rm
in}}\e\bbbr^3\times[0,\infty[ ,\tag 1.1$$ 
$$u\big|_{t=0}=u_0,\qquad
u_t\big|_{t=0}=u_1,\tag 1.2$$ 
where $g:\bbbr\to\bbbr$ and the initial data are
given sufficiently smooth functions, and 
$u_t={\partial\over \partial t}u$, etc. The linear case 
$g(u)=mu$, where $m\in
\bbbr$, corresponds to the classical Klein Gordon 
equation in relativistic particle
physics; the constant $m$ may be interpreted as a mass 
and hence is generally
assumed to be nonnegative. In an attempt to model also 
nonlinear phenomena 
like quantization, in the 1950s equations of type (1.1) 
with nonlinearities like
$$g(u)=mu+u^3,\qquad m\ge 0,$$
were proposed as models in 
relativistic
quantum mechanics with local interaction; see for 
instance Schiff~[13] 
and
Segal~[14]. Solutions could be real or complex-valued 
functions. In the latter case it
was natural to assume that the nonlinearity commutes with 
the phase, that is,
there holds
$$g(e^{i\va}u)=e^{i\va}g(u)\qquad\hbox{for 
all}\e\va\in\bbbr,$$
and hence, in particular, that
$g(0)=0$.
In this case, $g$ may be expressed
$$g(u)=u\th f\bigl(|u|^2\bigr),$$
which gives the form of equation (1.1) studied, for 
instance, by J\"orgens~[8]. Here,
for simplicity, and since all important features of our 
problem already seem to
exist in this case, we confine ourselves to the study of 
real-valued solutions of
equation~(1.1).
To model  effects thought to arise in the case, for 
instance, of spinor
fields $u$, the scalar equation~(1.1)  also  has been 
considered in space dimensions
$n\ge 3$; see~[14]. 

Various other models involving nonlinearities $g$ 
depending also on $u_t$ and
$\nabla u$, the spatial gradient of $u$, have been 
studied. The so-called
``$\sigma$-model" involves an equation of type (1.1) for 
vector-valued functions
subject to a certain (nonlinear) 
constraint.\footnote"$^1$"{In fact,
as observed by Shatah and Tahvildar-Zadeh \cite{21},
under suitable symmetry assumptions also
$\sigma$-models give rise to semilinear wave equations of 
type
(1.1) on $\Bbb{R}^4\times\Bbb{R}$.}
 In this case  $$g(u)=u\bigl(|u_t|^2-|\nabla
u|^2\bigr),$$ and the solution $u=(u_1,\dots,u_n)$ is 
constrained to satisfy
the condition $$|u|^2=u^2_1+\dots+u^2_n=1;$$
see Shatah~[15] for some recent results on this problem 
and references. 

To limit this
survey to a reasonable length, however, we restrict our 
study
to nonlinearities depending only on $u$; that is, the 
semi-linear case. The examples
stated previously suggest that we assume that $g(0)=0$ 
and that $g$ satisfies
polynomial growth 
$$
\bigl|g(u)\bigr|\le C(1+|u|^{p-2})|u|\quad\hbox{for some}\e
p\ge 2,\ C\in\bbbr.\tag1.3
$$ 
Moreover, following Strauss~[16, Theorem 3.1], we will 
assume that $g$ satisfies the
conditions
$$G(u)\ge -C|u|^2\quad\hbox{for some}\e C\in\bbbr\th,\tag 
1.4$$
and
$$\bigl|G(u)\bigr|{\di 
/}\bigl|g(u)\bigr|\RA\infty\quad\roman{
as}\e|u|\RA\infty\th,\tag 1.5$$
where
$G(u)=\int^u_0 g(v)dv$. Let us briefly motivate the 
latter two conditions.

First, (1.4) 
and (1.5) include the linear case (with no sign 
condition) or, more generally, the
case of Lipschitz nonlinearities. Second, in the 
super-linear case, that is, if
$\bigl|g(u)\bigr|{\di /}|u|\to \infty$ 
as $|u|\to\infty$, conditions (1.4), (1.5) should be 
regarded as a coerciveness
condition. In fact, in this case finite propagation speed 
$\le 1$ and conservation of
energy imply locally uniform a priori bounds in $L^2$ for 
solutions of (1.1) in terms
of the initial data; this will be developed in detail in 
\S 2. 

By contrast, in
the noncoercive case it is easy to construct solutions of 
(1.1) with smooth initial
data that blow up in finite time; for instance, for any 
$\a>0$ the function
$$u(x,t)={1\over (1-t)^\a}$$ solves the equation 
$$u_{tt}-\2
u=\a(1+\a)u|u|^{2\over\a}$$ and blows up at
$t=1$. Observe that for
$\a={1\over m}$, $m\in\Bbb N$, the right member of this 
equation is analytic.
Modifying the initial data off
$\bigl\{x;|x|\le 2\bigr\}$,
say, we even obtain a singular solution with
$C^\infty$-data having compact support. (See John~[7] for 
a blow-up result for a
similar equation.) 
Thus, conditions like (1.3)--(1.5) seem natural if we are 
interested in global solutions.

The class (1.3)--(1.5) includes the following special cases
$$g(u)=m u|u|^{q-2}+u|u|^{p-2},\qquad m\ge 0,\ 2\le
q<p\th.\tag 1.6$$ 

As we shall see, for nonlinearities of this kind the 
answer to the
existence problem for (1.1), (1.2) in a striking way 
depends on the space dimension
$n$ and on the exponent~$p$. In particular, in the 
physically interesting case $n=3$,
global existence for $p<6$ can be established with 
relative ease, while the same
question for $p>6$ so far has eluded all research 
attempts. The ``critical" case $p=6$
has only recently been settled 
and a comprehensive account of this result is one of the 
objectives pursued in this
survey.

In fact, the apparent existence of a ``critical power" 
for (1.1) and recent advances on
elliptic problems involving critical nonlinearities 
prompted our interest in the
$u^5$-Klein Gordon equation. ``Critical powers" very 
often come into play in
nonlinear problems through Sobolev embedding. In 
particular, $p=6$ is the critical
power for the Sobolev embedding $H^{1,2}
_{\operatorname{loc}}(\bbbr^3)\hookrightarrow
L^p_{\operatorname{loc}}(\bbbr^3)$. (In $n$ dimensions 
the critical power for this embedding is
$p={2n\over n-2}$.) Moreover, they very often arise 
naturally from the
requirements of scale invariance, that is, whenever 
``intrinsic" notions are involved.
A beautiful example of such a problem is the Yamabe 
problem concerning the
existence of conformal metrics with constant scalar 
curvature on a given (compact)
Riemannian manifold. Through the work of Trudinger, 
Aubin, and---finally---Schoen
this problem has now been completely solved and it has 
become apparent that at the
critical power properties like ``compactness of the 
solution set" 
depend crucially on
global aspects of the problem; in this case, on the 
topological and differentiable
structure of the manifold. See Lee and Parker [9] for a 
recent survey of the Yamabe
problem in this journal. 

Incidentally, for nonlinear wave
equations (or nonlinear Schr\"odinger equations $iu_t-\2 
u+u|u|^{p-2}=0$) there
appear to be many ``critical powers,'' depending on what 
aspect of the problem we
consider: global existence, scattering theory,  \dots  ; 
see Strauss~[16,~p.~14f.]. As
regards global existence, it remains to be seen whether 
the critical power represents
only a technical barrier or, in fact, defines the 
dividing line between qualitatively
different regimes of behavior of (1.1), (1.2). Through 
this survey I would like to
invite further research on this topic. 

We
conclude this introduction with a short overview of the 
existence results in the case
of a pure power $$u_{tt}-\2 u+u|u|^{p-2}=0,\qquad 
p>2.\tag 1.7$$ 

For more general
nonlinearities of type (1.6) similar results hold true. 
(In contrast, for problems
related to scattering, also the lower order terms of $g$ 
may be decisive.)  

\subheading{The sub-critical case}
For $n=3$,
$p<6$
global existence and regularity was established by 
J\"orgens~[8] in 1961. J\"orgens
also was able to show local (small time) existence of 
regular solutions to (1.7), (1.2)
for arbitrarily large
$p$.
Moreover, he was able to reduce the problem of existence 
of global, regular
solutions to (1.1) to (local) estimates of the
$L^\infty$-norms of solutions.
%

These results have been generalized to higher dimensions; 
however, such
extensions have been very hard to obtain.
While J\"orgens' work relies on the classical 
representation formula for the
3-di\-mensional wave equation, this method fails in 
higher dimensions
$n>3$.
The fundamental solution to the wave equation no longer 
is positive;
moreover, it carries derivatives transverse to the wave 
cone. Nevertheless,
at least for
$n\le 9$,
the existence results of Pecher~[11], Brenner-von 
Wahl~[2] now cover the full
sub-critical range
$p<{2n\over n-2}$.  Regular solutions are unique.
\subheading{Global weak solutions}
On the other hand, by a suitable
approximation  and using energy estimates,
for all $p>2$, $n\ge 3$
it is possible to construct global weak solutions, 
satisfying (1.7) in a distributional
sense; see Segal~[14],
Lions~[10]. In this case, it even suffices to assume that 
the initial data
$u_0,u_1\th \in L^2_{\operatorname{loc}}(\bbbr^n)$
with $u_0\in L^p_{\roman{loc}}(\Bbb{R}^n)$ and 
distributional derivative
$\nabla u_0\in L^2_{\operatorname{loc}}(\bbbr^n)$.
Energy estimates immediately give uniqueness of weak 
solutions in case
${p\le{2n\over n-2}}-{2\over n-2}$; see Browder~[3]. 
However, this range is well
below the critical Sobolev exponent
$p={2n\over n-2}$. In order to improve the range of 
admissable exponents, more
sophisticated tools were developed, based, in particular, 
on the
$L^p-L^q$-estimates for the wave operator by 
Strichartz~[17]; see also Brenner~[1].
In their simplest version, these estimates allow to prove 
uniqueness of solutions to
(1.7), (1.2) for $p\le{2(n+1)\over n-1}$, the Sobolev 
exponent in
$(n+1)$
space dimensions.  
In fact, uniqueness can be
established for $p<{2n\over n-2}$; see Ginibre-Velo~[4]. 
In this case, moreover, the
unique solution can be shown to be ``strong," that is, to 
possess second derivatives in
$L^2$
and to satisfy the energy identity [4]. Some of these 
results will be
derived in \S2. 
\subheading{The critical case}
In dimension
$n=3$, global existence of
$C^2$-solutions in the critical case
$p=6$
was first obtained by Rauch~[12], assuming the initial 
energy
$$E\bigl(u(0)\bigr)=\int_{\bbbr^3}\left({|u_1|^2+|\nabla 
u_0|^2\over 2}+{|u|^6\over
6}\right)dx$$
to be small. His results will be presented in \S3.

In 1987, also for ``large" data global $C^2$-solutions 
were shown to exist by this
author~[18] in the radially symmetric case
$u_0(x)=u_0\bigl(|x|\bigr)$, $u_1(x)=u_1\bigl(|x|\bigr)$. 
Finally, Grillakis~[6] in 1989
was able to remove the latter symmetry assumption, 
yielding the following result:
\thm{Theorem 1.1}
For any
$u_0\in C^3(\bbbr^3)$,
$u_1\in C^2(\bbbr^3)$
there exists a unique solution
$u\in C^2\bigl(\bbbr^3\times[0,\infty[\bigr)$
to the Cauchy problem
$$u_{tt}-\2 u+u^5=0,\tag 1.8$$
$${u}_{{\big|}t=0}=u_0,\qquad {u_t}_{{\big|}t=0}=u_1.\tag 
1.9$$
\ethm

\n
In \S4 we present the detailed proof.
Related partial regularity
 results independently  have been obtained by Kapitanskii
\cite{20} in 1989. Uniqueness holds among
$C^2$-solutions. The proof procedes via a priori 
estimates. The classical 
representation formula crucially enters.
It seems unlikely that regularity or uniqueness of weak
solutions to (1.8), (1.9) can be established in a similar 
way.
Research on the critical case in higher dimensions is in 
progress; however, to this
moment the results on this subject still seem incomplete. 
Advances in these
questions may require eliminating the 
use of the wave kernel. 
\subheading{The super-critical case}
In
\S5  we observe that for sufficiently small initial data 
the existence of global
regular solutions, for instance, to the equation 
$$u_{tt}-\2
u+u^5+u|u|^{p-2}=0\quad\hbox{in}\th 
\bbbr^3\times[0,\infty[,$$ 
for any
$p>2$
can be deduced as a corollary to Rauch's result. Various
qualitative properties of solutions in the super-critical 
case have recently been
studied by Zheng~[19].

Other open problems concern scattering theory, involving, 
in particular, decay
estimates for solutions of (1.1) (see Ginibre-Velo~[5]), 
or existence and regularity
results for initial-boundary value problems.
\heading2. Preliminaries\endheading 

We begin our study of (1.1) with some general comments 
about local solvability and
global continuation of solutions to (1.1), (1.2). An 
excellent reference for many
fundamental results on nonlinear wave equations is 
Strauss~[16]; our treatment of
these issues will be somewhat narrower and directed 
towards our final goal: the
critical power. This restricted aim, however, will enable 
us to make this paper
essentially self-contained and to present a lot of 
material connected with the
existence problem for (1.1), (1.2) in detail, introducing 
the reader to various
approaches to this problem and showing their strengths 
and limitations. 
\subheading{Representation formulas} 
The representation of solutions to the
inhomogeneous wave equation in terms of the fundamental 
solution and energy
estimates form the basis of our solution method. For any 
$f\in C^\infty$, $\th
u_0,u_1\in C^\infty$ there exists a unique 
$C^\infty$-solution to the Cauchy
problem 
$$u_{tt}-\2 
u=f\quad\hbox{in}\e\bbbr^n\times[0,\infty[,\tag 2.1$$
$${u}_{{\big|}t=0}=u_0,\qquad {u_t}_{{\big|}t=0}=u_1.\tag 
2.2$$ 

\n
If $n=3$, the most
interesting case, this solution, in fact, is given by 
$$
\split
u(x,t)=&{d\over
dt}\left({1\over 4\pi t}\int_{\partial 
B_t(x)}u_0(y)\,dy\right)+ {1\over 4\pi
t}\int_{\partial B_t(x)}u_1(y)\,dy\\
 &+{1\over4\pi}\int^t_0\int_{\partial
B_{t-s}(x)}{f(y,s)\over t-s}\,dy\, ds,\endsplit
\tag 2.3$$ 
where
$B_r(x)=\bigl\{ y\in\rn;|x-y|<r\bigr\}$. 
From (2.3) we see immediately that information propagates 
with speed
$\le 1$.
In particular,
$u(t)$
has compact support for any
$t\ge 0$
if this is the case for $u_0,u_1$,
and
$f$. However, (2.3) also shows a fundamental weakness of 
the classical
approach: For $u_0\in C^3$, $u_1\in C^2$, 
$f\in C^2$, the solution $u$ will lie in $C^2$, only.
That is, we encounter a loss of differentiability.  
In higher dimensions, a representation formula similar 
to~(2.3) holds, however,
involving an even larger number (the integer part of 
$n\over 2$, resp. ${n-2\over
2}$ ) of derivatives of $u_0$, resp. of $u_1$ and $f$. 
This makes the representation
formula appear to be ill-suited for proving existence of 
solutions for semilinear
equations in dimensions $n>3$. 

By contrast, no loss of differentiability will occur if
instead of pointwise control of the solution we are 
content with control of integral
norms. The basic observation is the following. 
\subheading{Energy inequality}
Upon multiplying (2.1) by $u_t$ we obtain 
$$
{d\over
dt}\left({|u_t|^2+|\nabla u|^2\over 2}\right)-
\operatorname{div}
(u_t\nabla u)=f\th u_t,$$ where the
terms $$e_0(u)={|u_t|^2+|\nabla u|^2\over 2}$$ and 
$p(u)=u_t\nabla u$ may be
interpreted as energy and momentum of the solution $u$. 
Integrating
in $x$, if $u(t)$ has compact support, by H\"older's 
inequality we obtain
$$\eqalign{{d\over 
dt}E_0\bigl(u(t)\bigr)&\le\left(\int_{\bbbr^n}\bigl|
f(\cdot,t)\bigr|^2\,dx\right)^{1/2}\left(\int_{\bbbr^n}%
\bigl|
u_t(\cdot,t)\bigr|^2\,dx\right)^{1/2}\cr &\le\Bigl(2
E_0\bigl(u(t)\bigr)\Bigr)^{1/2}\bigl|\bigr| 
f(\cdot,t)\bigl|\bigr|_{L^2(\Bbb
R^n)},\cr}$$ where
$$E_0\bigl( 
u(t)\bigr)=\int_{\bbbr^n}e_0\bigl(u(\cdot,t)\bigr)\,dx=:%
\bigl|\bigr|
u(t)\bigl|\bigr|^2_0$$
denotes the ``energy norm."

\n
Thus
$${d\over dt}\bigl|\bigr|
u(t)\bigl|\bigr|_0\le{1\over\sqrt
2}\bigl|\bigr|f(\cdot,t)\bigl|\bigr|_{L^2(\bbbr^n)}\le%
\bigl|\bigr|
f(\cdot,t)\bigl|\bigr|_{L^2(\bbbr^n)}.\tag2.4$$ 
In particular, if $f=0$, the ``energy" $E_0$ is conserved.

Various other conservation laws
can be obtained by using further multipliers related to 
symmetries of the wave
operator. Very subtle identities and integral estimates 
in this way have been
found; see Strauss~[16, Chapter 2] for an overview of 
results. In particular, in
\S4 we will make use of the integral estimate implied by 
invariance of the
wave operator under dilations $(x,t)\mapsto(R\th x,R\th 
t)$ for $R>0$. For our
immediate uses, however, the energy inequality will 
suffice.  

So far, (2.4) has been
established rigourously only for $C^\infty$-data $u_0$,
$u_1$, and $f$ with spatially
compact support. For our next topic it is essential to 
extend the validity of (2.4) to
distribution solutions of~(2.1) for finite energy initial 
data, that is,
 for $u_0$, $u_1\in
L^2(\bbbr^n)$ with $\nabla u_0\in L^2(\bbbr^n)$, and 
functions $f$ belonging to
$L^2\bigl(\bbbr^n\times[0,T]\bigr)$ for any $T>0$. 
To achieve this extension, by
density of $C^\infty_0(\bbbr^n)$ in $L^2(\bbbr^n)$ we may 
approximate data
$u_0,u_1$ as above by functions $u_0^{(m)}$,
$u_1^{(m)}\in C^\infty_0$, converging to
$u_0,u_1$ in energy norm as $m\to\infty$. Similarly, for 
any $T>0$ we can find
smooth, compactly supported functions $f^{(m)}$ 
converging to $f$ in
$L^2\bigl(\bbbr^n\times[0,T]\bigr)$. Let $u^{(m)}$ be the
corresponding sequence of solutions to (2.1), (2.2), 
given by the classical
representation formula. Then, applying~(2.4) to the 
difference $v=u^{(m)}-u^{(l)}$
of any two solutions, we see that
$u^{(m)}(\cdot,t)$ is a Cauchy sequence in energy norm, 
uniformly in $t\in[0,T]$. The
limit $u$ is a distribution solution to (2.1), (2.2) with 
uniformly finite energy
in the interval $[0,T]$, which satisfies (2.4) in the 
slightly weaker sense
$$\3 u(t)\3_0\le\3
u(0)\3_0+\int^t_0\3 f(\cdot,s)\3_{L^2(\bbbr^n)}\,ds,\tag 
2.5$$
for all $t\le T$. In particular, $u$
is unique in this class. 

In a similar way, we now use (2.5) to
construct solutions to nonlinear wave equations (1.1), 
(1.2) for smooth, compactly
supported initial data and smooth nonlinearities 
satisfying a global Lipschitz
condition by a contraction mapping argument. 
\subheading{Global solutions for Lipschitz nonlinearities}
Indeed, if $g\:\Bbb
R\to\bbbr$ is smooth and globally Lipschitz, for any 
$v\in C^\infty_0
\bigl(\Bbb
R^n\times[0,\infty[\bigr)$ we obtain a 
$C^\infty$-solution $u=K(v)$ to the initial
value problem $$u_{tt}-\2 u=-g(v)$$ with data $u_0,u_1$. 
By (2.5), for all $T>0$ we have
$$
\eqalign{\sup_{0\le
t\le T}\Bigl|\Bigr|\bigl(K(v)-K(\tilde
v)\bigr)(t){\Bigl|\Bigr|}_0&\le\int^T_0{\Bigl|\Bigr|}%
\bigl( 
g(v)-g(\tilde
v)\bigr)(t)\SS_{L^2}\,dt\cr &\le L\int^T_0\3(v-\tilde 
v)(t)\3_{L^2}\,dt\cr &\le T\th
L\sup_{0\le t\le T}\3(v-\tilde v)(t)\3_{L^2},\cr}$$ where 
$L$
denotes the Lipschitz constant of
$g$. Moreover, if
$u_0,u_1$
have support in
$B_R(0)$,
and if
$v(t)$
has support in
$B_{R+t}(0)$,
so will
$u(t)$.
Finally, by Poincar\' e's inequality, for such $v,\tilde 
v$, and
$t\le 1$
we can estimate
$$\3(v-\tilde v)(t)\3_{L^2}\le(R+t)\3\nabla(v-\tilde 
v)(t)\3_{L^2}\le\sqrt
2(R+1)\3(v-\tilde v)(t)\3_0.$$
Thus, for
$T\le \min\left\{1,{1\over \sqrt 2 L(R+1)}\right\}$
the map
$K$
extends to a contracting map on the space
$$\nomultlinegap\multline
V=\Bigl\{ v\in L^2\bigl(\bbbr^n\times[0,T]\bigr);
\operatorname{supp}
 \bigl(v(t)\bigr)\subset B_{R+t}(0),v_t(t),\nabla v(t)\in 
L^2(\bbbr^n)\\
\hbox{for almost every}\e t,\e
\hbox{and}\e \sup\limits_{0\le t\le T}\3 
v(t)\3_0<\infty\Bigr\},\endmultline
$$
endowed with the norm
$$||v||_{_{\di V}}=\sup_{0\le t\le T}\3 v(t)\3_0.$$
Let
$u$
be the unique fixed point of
$K$
in
$V$; then $u$
weakly solves~(1.1)
and assumes its initial data (1.2) in the distribution 
sense. By an approximation
argument as in the preceding paragraph, likewise for 
compactly supported
measurable initial data with finite 
 energy we obtain a (unique) local solution to (1.1), 
(1.2) in the space $V$. Observe
that the support of the solution grows with speed~$\le 
1$. Hence, given any number
$T_0>0$, by iterating the above construction a finite 
number of 
times (with $T\le\min
\left\{1,{1\over \sqrt 2 L(R+T_0)}\right\}$ ) we obtain a 
finite energy solution to
(1.1), (1.2) on the interval $[0,T_0]$ for any finite 
energy initial data supported in a
ball of radius~$R$. Since $T_0$ is arbitrary, this 
solution can be continued globally.

Finally, by finiteness of propagation speed also the 
assumption that the initial data be
compactly supported can be removed. Indeed, if
$u$ and $v$ solve~(1.1) for  compactly supported, finite 
energy initial data that
coincide on the ball $B_R(0)$, their difference~$\tilde 
u$ will solve an equation of
type~(1.1) with a Lipschitz nonlinearity $\tilde g$, 
where $\tilde g(\tilde
u)=g(u)-g(v)$, and initial data vanishing on $B_R(0)$. By 
the above existence and
uniqueness result, $\tilde u$ is supported outside the 
light cone
$K=\bigl\{(x,t);|x|<R-t\bigr\}$ above $B_R(0)$. Hence a 
solution to~(1.1) on~$K$ is
entirely determined by its data in $B_R(0)$. For 
arbitrary data $u_0,u_1$ with
locally finite energy, and $k\in\bbbn$, we then let 
$u^{(k)}_0,u_1^{(k)}$ be compactly
supported data that agree with $u_0,u_1$ on
$B_k(0)$. For any $k\in\bbbn$ the corresponding global 
solutions $u^{(k)},k\ge k_0$,
then agree on the cone above $B_{k_0}(0)$. Hence the 
sequence $(u^{(k)})$
converges locally in energy norm to a global solution $u$ 
of (1.1), (1.2). 
In the same way, as far as global existence is concerned, 
in the
following for convenience---and with no loss of 
generality---we
may suppose that the
initial data have compact support. Moreover, for our next 
topics
 we also require the
data $u_0,u_1$ to be smooth. 
\subheading{Strong solutions}
Taking
difference quotients\th $u^{(h)}={u(\cdot)-u(\cdot+
he)\over h}$\th
in any space direction
$e$
and passing to the limit
$h\to 0$
we see that
$v=e\cdot \nabla u$
weakly solves
$$v_{tt}-\2 v+ g'(u)v=0$$
and satisfies
$$
\split
\3 v(T)\3_0-\3
v(0)||_0\le&\int^T_0\|\bigl(g'(u)v\bigr)(t)\|_{L^2}\,dt\le 
L\negthinspace
\int^T_0
\negthinspace\negthinspace \3
v(t)\3_{L^2}\,dt\\
\le& C\th L\il^T_0\3 v(t)\3_0\, dt,
\endsplit
\tag2.6
$$ 
for any $T>0$.
Thus
$\nabla u_t(t),\nabla^2 u(t)\in L^2(\bbbr^n)$,
uniformly in
$t\in[0,T]$, and from equation (1.1) it now also follows 
that
$u_{tt}(t)\in L^2$,
uniformly in
$t\in[0,T]$, for any $T>0$. 
This is the class of ``strong solutions" to nonlinear 
wave equations. 
For strong
solutions we can derive the strong form of the energy 
inequality~(2.4).
Since $g(u)u_t={d\over dt}G(u)$, upon multiplying~(1.1) 
by $u_t$ we obtain the
conservation law 
$${d\over dt}\left({|u_t|^2+|\nabla u|^2\over 2}+
G(u)\right) - 
\operatorname{div}(u_t\nabla
u)=0,$$
where the term
$$e\bigl(u(t)\bigr)={|u_t|^2+|\nabla u|^2\over 2}+
G(u)=e_0(u)+G(u)$$
now also contains the ``potential energy density" $G(u)$. 
Let
$$ E\bigl(u(t)\bigr)=\int_{\bbbr^n}e\bigl(u(t)\bigr)\,dx.$$
Integrating over $\bbbr^n$, since~$u(t)$ has compact 
support, we thus obtain
that
$${d\over dt}E\bigl(u(t)\bigr)=0,\tag 2.7$$
and energy is strictly conserved.
\subheading{Higher regularity}
By iterating the above procedure, we may
want to derive $L^2$-bounds for higher and
higher derivatives $D^ku$, where
$D$ denotes any space-time derivative, $k\in \Bbb N_0$.
For instance, in case
$k=3$, any second order spatial derivative $w=\nabla^2u$
satisfies
$$w_{tt}-\2 w+ g'(u)w+g''(u)|\nabla u|^2=0.$$
However, while the  term involving $g'(u)w$
 can be dealt with as before, the second
term presents some difficulty and can only be controlled 
in terms of
$E_0(w)$
if the dimension 
$n\le 8$. In this case, assuming $|g''(u)|\le C$,
by Sobolev's inequality we can estimate
$$
\SS g''(u)\bigl|\nabla u(t)\bigr|^2\SS_{L^2}\le
 C\3\nabla u(t)\3^2_{L^4}\le
C\3\nabla^2 u(t)\3^{2\gamma}_{L^2}\3\nabla^3 
u(t)\3^{2-2\gamma}_{L^2}
\le
C\th E_0(w)^{1-\gamma}
$$
where
${\gamma\over {2n\over n-2}}+{(1-\gamma)\over{2n\over 
n-4}}\le{1\over 4}$.
The energy inequality~(2.4) formally yields
$${d\over dt}\3 w(t)\3_0\le C\3 w(t)\3_0+C\3 
w(t)\3_0^{2-2\gamma}.$$
Note that
$\gamma\to 0$ as $n\nearrow 8$. Hence the last exponent 
may be $>1$, and
$\3 w(t)\3_0$ might blow up in finite time. Similar 
problems arise if we want to
control higher derivatives of~$u$ by this simple trick.

If
$n\ge 9$, we cannot start our iteration at~$k_0=2$. 
However, if we choose
$k_0\in\bbbn$ sufficiently large, by using Sobolev's 
embedding theorem as above
we can derive a priori estimates for $\3 D^k u(t)\3_0$ 
for any $k\ge k_0$ in a small
time interval $0\le t\le T(k)$. As in low dimensions, 
these estimates also may blow
up in finite time.

Nevertheless, we can use these estimates to show the 
local (small time) existence of
solutions to general semilinear equations~(1.1) for 
smooth, compactly supported
initial data.
\subheading{Local solutions for semi-linear equations}
Indeed, given an arbitrary smooth map
$g:\bbbr\to\bbbr$
we may approximate $g$
by maps
$g^{(k)}$
satisfying a uniform Lipschitz condition and coinciding 
with
$g$
for
$|u|\le k$.

By the preceding discussion, given any smooth initial 
data of compact support, for
each
$k\in\Bbb N$
we obtain a global strong solution
$u^{(k)}$
of the approximate equation
$$u_{tt}-\2 u+ g^{(k)}(u)=0,$$
with $\3 D^l u^{(k)}(t)\3_0\le C$
for any
$l\in\{0,\dots,l_0\}$
on some interval $0\le t\le T=T(l_0)$, where
$C$
depends on the Lipschitz constant of
$g^{(k)}$, $l_0$, $T$, and the size of the support of the 
initial data. 

If
$n=3$, by the Sobolev embedding theorem
$$\3 D^l u^{(k)}(t)\3_{L^\infty}\le C\3 
D^{2+l}u^{(k)}(t)\3_{L^2}\le C\3
D^{1+l}u^{(k)}(t)\3_0,$$
for $l=0,1,2$. In particular, for large $k\in\Bbb N$
and sufficiently small
$T>0$, we obtain
$\bigl|u^{(k)}(x,t)\bigr|\le k$
in
$\bbbr^3\times[0,T]$,
and $u^{(k)}$
will be a solution to (1.1). 
Similarly, if
$n>3$
we can bound
$$\3 D^l u^{(k)}(t)\3_{L^\infty}$$
in terms of 
$\3 D^{m+l}u^{(«k)}(t)||_0$, where
$m>{n\over 2}-1$, and we may conclude as before.

Again remark that by finiteness of propagation speed the 
assumption that the initial
data be compactly supported can be removed; in this case, 
however, we can only
assert the existence of a solution to (1.1), (1.2) in a 
neighborhood of
$\bbbr^n\times\{0\}$. 
\subheading{Global weak solutions}
We now specialize our
nonlinearity $g$ to be of the form (1.3)--(1.5). By 
assumption (1.4) there exist
sequences $r^\pm_k\RA\pm\infty$ as $k\to\infty$ such that 
$$r^\pm_kg\left(r^\pm_k\right)\ge 
-C\left|r^\pm_k\right|^2\th.$$
We approximate $g$ by Lipschitz functions
$$g^{(k)}(u)=\cases
g(r^-_k),&\text{if }u<r^-_k,\\
g(u),&\text{if }r^-_k\le u\le r^+_k,\\
g(r^+_k),&\text{if }u>r^+_k,\endcases
$$
with primitive $G^{(k)}(u)$. Note that the approximating 
functions $g^{(k)}$
satisfy~(1.4) with a uniform constant~$C$. Now, for
any $k\in\Bbb N$ and
smooth, compactly supported data we obtain a unique 
global strong solution $u^{(k)}$
to the approximate problem (1.1), (1.2) with
$D^2u^{(k)}(t)\in L^2(\bbbr^n)$ for all $t$.

Conservation of energy~(2.7) implies uniform bounds for 
$u=u^{(k)}$. Let
$$E^{(k)}\bigl(u(t)\bigr)=E_0\bigl(u(t)\bigr)+
\int_{\bbbr^n}G^{(k)}\bigl(u(t)\bigr)\,dx.$$
By (1.4), (2.7), for any $t\ge 0$ we have
$$E_0\bigl(u(t)\bigr)-C\3 u(t)\3^2_{L^2}\le
E^{(k)}\bigl(u(t)\bigr)=E^{(k)}\bigl(u(0)\bigr)\le 
C<\infty,\tag 2.8$$
uniformly in $k\in \bbbn$. In order to control $\3 
u(t)\3_{L^2}$, for 
fixed $x\in\bbbr^n$ we estimate
$$
\bigl|u(x,t)-u_0(x)\bigr|^2=\left|\int^t_0u_t(x,s)\,ds%
\right|^2\le
t\int^t_0\bigl|u_t(x,s)\bigr|^2\,ds.$$
Integrating in $x$, by Minkowski's inequality we obtain
$$\3 u(t)\3_{L^2}\le ||u_0||_{L^2}+\left(2t\int^t_0 
E_0\bigl(u(s)\bigr)\,ds\right)^{1/2}.$$

For $t\le T$
this and~(2.8) gives the integral inequality
$$E_0\bigl(u(t)\bigr)\le C+CT\int^t_0 
E_0\bigl(u(s)\bigr)\,ds$$
for
$E_0\bigl(u(t)\bigr)$. From Gronwall's lemma we thus 
conclude that $u(t)=u^{(k)}(t)$
is uniformly bounded in energy norm on any interval 
$[0,T]$, uniformly in
$k\in\bbbn$. 

Hence,
$(u^{(k)})_{k\in\Bbb N}$
is weakly relatively compact in the energy norm. 
Moreover, the support of $u^{(k)}(t)$
is bounded uniformly in
$k$, for all
$t\le T$.
By the Rellich-Kondrakov
theorem, therefore, we may assume that
$u^{(k)}\to u$
strongly in $L^2(Q)$
on any compact space-time region $Q$ and pointwise almost 
everywhere. The limit
$u$ has finite energy
$$E_0\bigl(u(t)\bigr)\le\liminf_{k\to\infty}E_0%
\bigl(u^{(k)}(t)\bigr),$$
and
$$\int_{\bbbr^n}G\bigl(u(t)\bigr)\,dx\le\liminf_{k\to 
\infty}\int
_{\bbbr^n}G^{(k)}\bigl(
u^{(k)}(t)\bigr)\,dx\tag 2.9$$
for almost every
$t>0$,
by Fatou's lemma. Finally, for
$\va\in C^\infty_0\bigl(\bbbr^n\times]0,\infty[\bigr)$ we 
obtain
$$\eqalign{\int\int(u_t\va_t-\nabla u\th\nabla\va)\,dx\,
dt&=\lim_{k\to\infty}\int\int(u^{(k)}_t\va_t-\nabla 
u^{(k)}\nabla\va)\,dx\,
dt\cr
&=\lim_{k\to\infty}\int\int g^{(k)}(u^{(k)})\va\th dx\ dt
=\int\int g(u)\va\, dx\, dt,\cr}$$
where $\int\int\dots$ denotes integration over
$\bbbr^n\times[0,T]$. That is, $u$ weakly solves equation
(1.1). (Vitali's theorem, (1.5) and (2.9) were used to 
pass to the limit in the
nonlinear term.) Similarly, approximating
$L^2$-data
$u_0,u_1$
of finite energy
$$\int_{\bbbr^n}\left({|u_1|^2+|\nabla u_0|^2\over 2}+
G(u_0)\right)dx<\infty$$
by functions
$u_0^{(k)},u_1^{(k)}\in C_0^\infty(\bbbr^n)$,
the existence of global weak solutions to (1.1) for 
arbitrary finite energy data may be
derived.
\subheading{Regularity and uniqueness}
In the special case (1.7) with
$p\le{2n\over n-2}-{2\over n-2}$ energy estimates may be 
used to obtain higher
regularity and uniqueness. Indeed, let
$u^{(h)}$
be a difference quotient in direction $e\in\bbbr^n$
as before. Then, upon passing to the limit
$h\searrow 0$, for $v=e\cdot\nabla u$ we obtain
$$v_{tt}-\2 v=(1-p)|u|^{p-2}v,$$
and thus, formally, by H\"older's inequality and (2.4), 
that
$$
{d\over dt}\3 v(t)\|_0\le 
C\3\bigl(|u|^{p-2}v\bigr)(t)\|_{L^2}\le C\3
u(t)\3^{p-2}_{L^{2^*}}\3 v(t)\3_{L^{2^*}},$$ 
where
$2^*={2n\over n-2}$. Sobolev's inequality now implies that
$${d\over dt}\3 v(t)\3_0\le C\3 u(t)\3_0^{p-2}\3 v(t)\3_0$$
and it follows that
$E_0\bigl(v(t)\bigr)<\infty$
for all $t$; that is,
$u$ is a strong solution to~(1.1). Similarly, higher 
regularity (for small time, if the
dimension is large) may be obtained. To see uniqueness, 
let  $u,\tilde u$ be
solutions to (1.7) with the same initial data (1.2). For 
$v=u-\tilde u$ we obtain the
inequality $${d\over dt}\3 v(t)\3_0\le C\Bigl(\3 u(t)\3_0+
\3\tilde
u(t)\3_0\Bigr)^{p-2}\3 v(t)\3_0.$$
Since
$v(0)=v_t(0)=0$,
uniqueness follows. 

By more sophisticated methods the above regularity and
uniqueness results may be extended to the full 
sub-critical range
$p<{2n\over n-2}$; see Ginibre-Velo~[4]. One such method 
will be briefly
explained next.
\subheading{$L^p-L^q$-estimates}
By a result of
Strichartz [17], for any $L^p$-solution of the wave 
equation in
$\rn\times[0,\infty[$
with
$\square u=u_{tt}-\2 u\in L^{p\over p-1}$
and vanishing initial data we have
$$\3 u(t)\3_{L^p}\le C\int^t_0(t-s)^{1-2n\delta}\3\square 
u(s)\3_{L^{p\over
p-1}}\,ds,\tag 2.10$$
provided
$\delta={1\over 2}-{1\over p}\le{1\over n+1}$; see also 
Brenner~[1]. We
illustrate how this estimate may be used to obtain global 
strong (or even classical)
solutions to (1.7) in dimensions $n>3$.
Since we will need $g\in C^2$, we suppose that $p\ge 3$. 
The above condition on $\dt$
then requires $n\le 5$, $p\le{2(n+1)\over n-1}$.

It suffices to show existence of a solution on
$\rn\times[0,T]$
for compactly supported, smooth data and for arbitrary
$T>0$. Let $u^{(0)}$ solve
$u^{(0)}_{tt}-\2 u^{(0)}=0$
with initial data (1.2), and for $k\in\Bbb N$
let
$u^{(k)}$
be the solution to the approximate equation
$$u^{(k)}_{tt}-\2 u^{(k)}+u^{(k)}\min\bigl\{ 
|u^{(k)}|^{p-2},k^{p-2}\bigr\}=0,$$
together with (1.2). Here, $\min \{a,b\}$ is a smooth 
function coinciding with the
minimum of $a$ and $b$ for $|a-b|\ge 1$. 

Then
$u^{(k)}=u^{(0)}+v^{(k)}$,
where
$$\square v^{(k)}=v^{(k)}_{tt}-\2 
v^{(k)}=-u^{(k)}\min\bigl\{
|u^{(k)}|^{p-2},k^{p-2}\bigr\},\tag 2.11$$ with
$$v^{(k)}(0)=v^{(k)}_t(0)=0.$$
By (2.6),
$D^2u^{(k)}(t)\in L^2(\rn)$
for all
$t$.
Moreover, by (2.9) we can uniformly bound
$$E_0\bigl(u^{(k)}(t)\bigr)\le C_0,$$
$$\int_{\rn}\min\left\{{|u^{(k)}|^p\over p},{k^p\over
p}+{\bigl(|u^{(k)}|^2-k^2\bigr)k^{p-2}\over 
2}\right\}dx\le C_0.$$
That is,
$\square v^{(k)}(t)$
is uniformly bounded in
$L^{p\over p-1}(\rn)$
for all
$t$.
From (2.10) we now obtain
$$\eqalign{\3 u^{(k)}(t)\3_{L^p}&\le\3 u^{(0)}(t)\3_{L^p}+
\3 v^{(k)}(t)\3_{L^p}\cr
&\le C(t)+C\int^t_0(t-s)^{1-2n\dt}\3\square 
v^{(k)}(s)\3_{L^{p\over p-1}}\,dx
\le C(T)\cr}\tag 2.12$$
for all $t<T$, since
$\dt={1\over 2}-{1\over p}\le{1\over n+1}<{1\over n}$ for 
the range of
$p$ and $n$
considered.

Differentiating (2.11), similarly we obtain
$$\bigl|\square(D
v^{(k)})\bigr|\le(p-1)|Du^{(k)}|\min\bigl%
\{|u^{(k)}|^{p-2},k^{p-2}\bigr\}
\le(p-1)|D u^{(k)}|\th|u^{(k)}|^{p-2}.$$

Hence
{$\square (Dv^{(k)})\in L^{p\over p-1}$}
and by (2.12) 
$$\3\square(D v^{(k)})(t)\3_{L^{p\over p-1}}\le (p-1)\3
Du^{(k)}(t)\3_{L^p}\3 u^{(k)}(t)\3^{p-2}_{L^p}
\le C(T)\3 Du^{(k)}(t)\3_{L^p}.$$
Thus, from (2.10) we infer
$$\3 D u^{(k)}(t)\3_{L^p}\le C(T)+
C(T)\int^t_0(t-s)^{1-2n\dt}\3 D
u^{(k)}(s)\3_{L^p}\,ds,$$
and it follows that
$\3 D u^{(k)}(t)\3_{L^p}\le C(T)$
for all
$0\le t\le T$.

Finally, we have
$$\bigl|\square(D^2v^{(k)})\bigr|\le C\bigl(|D
u^{(k)}|^2|u^{(k)}|^{p-3}+
|D^2u^{(k)}||u^{(k)}|^{p-2}\bigr),$$
whence
$$\eqalign{\3\square(D^2v^{(k)})(t)\3_{L^{p\over p-1}}
&\le C\3 D
u^{(k)}(t)\3^2_{L^p}\3 u^{(k)}(t)\3^{p-3}_{L^p}+C\3 
D^2u^{(k)}(t)\3_{L^p}\3
u^{(k)}(t)\3^{p-2}_{L^p}\qquad\cr
&\le C(T)\bigl(1+\3 D^2u^{(k)}(t)\3_{L^p}\bigr),\cr}$$
and from (2.10) again it follows that
$\3 D^2u(k)(t)\3_{L^p}\le C(T)$,
uniformly in
$k$.
But by Sobolev's inequality, for
$n\le 5$, $p\ge 3$,
we may estimate
$$\3 u^{(k)}(t)\3_{L^\infty}\le C\3 
D^2u^{(k)}(t)\3_{L^p}\le C(T).$$
It follows that for sufficiently large $k$ the function 
$u=u^{(k)}$
is a (strong) solution to the original equation (1.7).
If $g\in C^4$, we can proceed to bound the first and 
second derivatives of $u$ and
hence obtain a classical solution. 

Note that the range $p\le{2n+2\over n-1}$,
where Stichartz' estimate may be applied, slightly 
exceeds the range $p\le{2n\over
n-2}-{2\over n-2}$,
where simple energy estimates suffice to show regularity 
and uniqueness.
\subheading{Classical solutions}
If $n=3$, using (2.3) one can also devise a
contraction mapping argument in the space $C^2$ to obtain 
local classical solutions
to (1.1), (1.2) for initial data $u_0\in C^3, u_1\in C^2$ 
with compact support. 

Indeed, via~(2.3) the initial value problem (1.1), (1.2) 
can be converted into the
integral equation
$$u(x,t)=v(x,t)-{1\over 4\pi}\int
^t_0\int_{|x-y|=t-s}{g\big(u(y,s)\big)\over t-s}\,dy\,
ds,$$ 
where
$v$ denotes the solution to the homogeneous wave equation 
with data~(1.2). If $g$
is smooth and globally Lipschitz this can easily be 
solved on
$\bbbr^3\times[0,T]$ for suitably small $T>0$ by a 
contraction mapping argument in
the space $C^0\big(\bbbr^3\times[0,T]\big)$
with the $L^\infty$-norm. Differentiating (1.1) in any 
spatial
direction, similarly we obtain 
$$
Du(x,t)=Dv(x,t)-{1\over 
4\pi}\int^t_0\int_{|x-y|=t-s}{g'(u)Du\over t-s}\,dy\,ds$$
and an analogous equation for the second 
spatial derivatives, from which we can as usual
derive locally uniform bounds for all first and second 
derivatives of $u$ on
$\bbbr^3\times[0,T]$. To extend $u$ beyond $T$ we 
write
$$u(x,t)=v_1(x,t)-{1\over4\pi}\int^t_T\int_{|x-y|=t-s}{g%
\big(u(y,s)\big)\over 
t-s}\,dy\,
ds,$$ where now
$v_1$ denotes the solution of the homogeneous wave 
equation with data
$u(\cdot,T)$ and $u_t(\cdot,T)$ at time $T$. At first it 
may seem as if we had lost one
derivative in this procedure. However, following 
J\"orgens~[8,~p.~301], we can write
$$v_1(x,t)=v(x,t)-{1\over 4\pi}\int^T_0\int
_{|x-y|=t-s}{g\big(u(y,s)\big)\over t-s}\,dy\,
ds,$$
and
$v_1\in C^2$, as desired. Thus, for smooth Lipschitz 
nonlinearities by iteration we
obtain global $C^2$-solutions. Likewise, for smooth $g$ 
we obtain local
$C^2$-solutions (for small time). However, if $g'(u)$ is 
uniformly bounded (for
instance, if $u$ is uniformly bounded) on any interval 
$[0,T]$, then also this solution
extends globally. 
Finally, by finiteness of propagation speed, the 
assumption that the data have
compact support can be dropped. 

Due to loss of differentiability in the nonlinear term, 
in dimensions
$n>3$ this approach---apparently---fails.

After this preliminary discussion
of different approaches to semi-linear wave equations we 
now focus our attention
on~(1.7) in the critical case $p=6$ in dimension $n=3$, 
which will be fixed from now
on.  

In the next section  we present the existence result of 
Rauch for small
data. Then we present an energy decay estimate and show 
how regularity in the
radial case may be derived. Finally we focus on the work 
of Grillakis~[6], whose
penetrating analysis provides the crucial insight needed 
to pass from the radially
symmetric to the general case and give a slightly 
simplified exposition. We conclude
this paper with a global existence result for certain 
super-critical nonlinearities and
small data.
\heading 3. Rauch's result\endheading
Let
$z=(x,t)$
denote points in space-time
$\bbbr^3\times\bbbr$.
Given
$z_0=(x_0,t_0)$,
let
$$K(z_0)=\Bigl\{ z=(x,t);\e|x-x_0|\le t_0-t\Bigr\}$$
be the backward light cone with vertex
at
$z_0$,
$$M(z_0)=\Bigl\{ z=(x,t);\e|x-x_0|= t_0-t\Bigr\}$$
its mantle,
$$D(t;z_0)=\Bigl\{ z=(x,t)\in 
K(z_0)\Bigr\}\qquad(t\e\roman{fixed})$$
its space-like sections. If $z_0=(0,0), z_0$
will be omitted. For any space-time region
$Q\subset\bbbr^3\times\bbbr$,
$T<S$,
we denote
$$Q^S_T=\Bigl\{z=(x,t)\in Q;T\le t\le S\Bigr\}$$
the trunctated region. Hence, for instance,
$$\partial K^s_t=D(s)\cup D(t)\cup M^s_t.$$
If
$s=\infty$
or
$t=-\infty$,
they will be omitted. Given a function
$u$
on a cone
$K(z_0)$,
we denote its energy density by
$$
e(u)=\tfrac12
\bigl(|u_t|^2+|\nabla u|^2\bigr)+\tfrac16
|u|^6,$$
and by
$$E\bigl(u;D(t;z_0)\bigr)=\int_{D(t;z_0)}e(u)\,dx$$
its energy on the space-like section
$D(t;z_0)$.
Moreover, let
$x=y+x_0$
and denote
$$
d_{z_0}(u)=\tfrac12
\left| {y\over|y|}u_t-\nabla u\right|^2+\tfrac16
|u|^6$$
the energy density of 
$u$
tangent to
$M(z_0)$.

Finally, for
$x_0\in\bbbr^3$
let
$$B_R(x_0)=\bigl\{x\in\bbbr^3;\e |x-x_0|<R\bigr\}$$
with boundary
$$S_R(x_0)=\bigl\{x\in\bbbr^3;\e |x-x_0|=R\bigr\}.$$
In the following, the letters
$c,C$
will denote various constants.
$E_0$
will denote a bound for the initial energy. 

The proof of Rauch's existence result
relies on the following inequalities of Hardy-type that 
also play an essential role in
the work of Grillakis and this author on the limit case
$p=6$. We state these estimates in a form due to 
Grillakis [6, Lemma 2.1]).
\thm{Lemma 3.1}
Suppose 
$u\in L^6(B_R)$
possesses a weak radial derivative
$u_r={x\cdot\nabla u\over |x|}\in L^2(B_R)$.
Then with an absolute constant
$C_0$
for all
$0\le\rho <R$
the following holds\RM:
$$\int_{B_R\backslash B_\rho}{|u(x)|^2\over |x|^2}\,dx\le
4\int_{B_R\backslash
B_\rho}|u_r|^2\,dx+2R^{-1}
\int_{S_R}|u|^2\,do;\tag i$$
$$\int_{B_R}{|u(x)|^2\over |x|^2}\,dx\le
C_0\left(\int_{B_R}|u_r|^2\,dx+
\left(\int_{B_R}|u|^6\,dx\right)^{1/3}\right);\tag ii$$
$$
\aligned
&\int_{S_R}|u|^4\,do\le
C_0\left(\left(\int_{B_R}|u_r|^2dx\right)^{1/2}
\left(\int_{B_R}|u|^6\,dx\right)^{1/2}
+\left(\int_{B_R}|u|^6\,dx\right)^{2/3}\right).
\endaligned
\tag iii$$
\ethm
\demo{Proof}
(i) follows upon integrating the inequality
$$|u_r|^2=\left|{1\over\sqrt r}\partial_r(\sqrt r 
u)-{1\over
2r}u\right|^2\ge{|u|^2\over 4r^2}-{1\over 
2r^2}{\partial\over \partial r}(r\th
u^2)$$
over
$B_R\backslash B_\rho$.
See Grillakis [6, Lemma 2.1] for the remaining details of 
the proof.
\qed\enddemo

Let
$z_0=(x_0,t_0)$
be given and suppose
$u$
is a
$C^2$-solution of (1.8) on the backward light cone 
$K_0(z_0)$.
As observed in \S2 above, in order to prove that
$u$
can be extended to a global solution of
(1.8), it suffices to show that for any
$z_0$
as above
$$m_0=\sup_{K_0(z_0)}|u|$$
can be a  priori bounded in terms of $z_0$ and the 
initial data.
Clearly, we may assume that
$m_0=\big|u(z_0)\big|$.

The first fundamental estimate towards deriving a priori 
bounds of this kind is the
following local version of the energy inequality. For 
later use we refer to a slightly
more general situation than described above. 
\thm{Lemma 3.2}
Let $\bar z=(\bar x,\bar t)$. Suppose
$u\in C^2\bigl(K_0(\bar z)\backslash\{\bar z\}\bigr)$
solves {\rm(1.8), (1.9)}. Then for any
$0\le t<s<\bar t$
there holds
$$
E\bigl(u;D(s;\bar z)\bigr)+{1\over\sqrt 
2}\int_{M_t^s(\bar z)}d_{\bar z}(u)\,do= E
(u;D(t;\bar z)\bigr)\le E_0,$$
where $do$ denotes the surface measure on $M(\bar z)$.
\ethm
\demo{Proof}
Integrate the identity
$$\bigl(u_{tt}-\2 u+u^5\bigr)u_t={d\over 
dt}e(u)-\operatorname{div}
(\nabla u\th u_t)=0$$
over
$K^s_t(\bar z)$.
Now let
$y=x-\bar x$
and use the fact that the outward unit normal on $M(\bar 
z)$ is given by
$$
n={1\over \sqrt 2}\left({y\over |y|},1\right)\th.$$
Hence the ``energy flux" through $M(\bar z)$ is given by
$$\align n\cdot\big(-\nabla u\th u_t,e(u)\big)
&={1\over\sqrt
2}\Biggr({1\over2}\left(|\partial_tu|^2-2{y\over 
|y|}\cdot\nabla 
u\th u_t+|\nabla u|^2\right)+{1\over 6}|u|^6\Biggr)\\
&={1\over\sqrt{2}}d_{\bar z}(u).\endalign$$
See Rauch [12].\qed\enddemo

By Lemma 3.2, for any fixed
$\bar z$
the energy
$E\bigl(u;D(s;\bar z)\bigr)$
is a monotone nonincreasing function of
$s\in[0,\bar t[$
and hence converges to a limit as
$s\nearrow\bar t$.
It follows that
$$
\int_{M^s_t(\bar z)}\bigl(d_{\bar z}(u)\bigr)do\to 
0\qquad (s,t\nearrow\bar
t);\tag 3.1$$
however, at a rate that may depend on
$\bar z$.

Following Rauch [12] we now decompose
$u=v+w$,
where
$v\in C^2\bigl(\bbbr^3\times[0,\infty[\bigr)$
is the unique solution of the homogeneous wave equation
$$v_{tt}-\2 v=0$$
with initial data (1.9) and
$$w_{tt}-\2 w+u^5=0,\qquad w_{{\di |}t=0}=0={w_t}_{{\di 
|}t=0}.$$
In particular, at
$z_0=(x_0,t_0)$
we may express 
$w$
via (2.3) as follows
$$w(z_0)=-{1\over 4\pi}\int_{M_0(z_0)}{u^5(x,t)\over 
t_0-t}\,do(x,t).$$

Thus, and splitting integration over
$M^T_0(z_0)$
and
$M_T(z_0)$
for suitable
$T$,
we obtain
$$\eqalign{m_0
&=\bigl| u(z_0)\bigr|\le\bigl| v(z_0)\bigr|+\bigl| 
w(z_0)\bigr|\cr
&\le \big| v(z_0)\big|+{m_0\over 4\pi}\int_{M_T(z_0)}
 {u^4\over t_0-t}\,do+{1\over 4\pi}
\int_{M^T_0(z_0)}
{|u|^5\over t_0-t}\,do.\cr}\tag 3.2$$ 

By H\"older's inequality and Lemma~3.2, the last term
$$\eqalign{\int_{M_0^T(z_0)}{|u|^5\over t_0-t}\,do&\le
C\big(|t_0-T|^{-1/2}-|t_0|^{-1/2}\big)\left(%
\int_{M_0(z_0)}|u|^6\,do\right)^{5/6}\cr 
& \le C\th E_0^{5/6}|t_0-T|^{-1/2}\th.\cr}\tag 3.3$$

Hence, if for some $T<t_0$ we have
$$\int_{M_T(z_0)}{u^4\over t_0-t}\,do\le 2\pi,\tag 3.4$$
from (3.2) and (3.3) we can bound
$$m_0\le 2\big|v(z_0)\big|+C\th E_0^{5/6}(|t_0-T|^{-1/2}
-|t_0|^{-1/2})\tag 3.5$$
in terms of the initial data, $z_0$, and $T$.

Now, by H\"older's inequality
$$\int_{M_T(z_0)}{u^4\over 
t_0-t}\,do\le\left(\int_{M_T(z_0)}{|u|^2\over
|t_0-t|^2}\,do\right)^{1/2}\left(\int_{M_T(z_0)}|u|^6\,do%
\right)^{1/2}.$$
Let $\tilde u(y)=u\bigl(x_0+y,t_0-|y|\bigr)$. Then by 
Lemma 3.1 we have
$$
\split
\int_{M_T(z_0)}{|u|^2\over
|t_0-t|^2}\,do&=\sqrt 2 \int
_{B_{t_0-T}(0)}{\bigl|\tilde{u}(y)\bigr|^2\over|y|^2}\,dy\\
&\le C\int_{B_{t_0-T}(0)}|\nabla \tilde u|^2\,dy+
C\left(\int
_{B_{t_0-T}(0)}|\tilde u|^6\,dy\right)^{1/3}\\
&\le
C\int_{M_T(z_0)}d_{z_0}(u)\,do+
C\left(\int_{M_T(z_0)}d_{z_0}(u)\,do\right)^{1/3}\\
&\le
C\Bigl(E\big(u;D(T;z_0)\big)+
E^{1/3}\big(u;D(T;z_0)\big)\Bigr).\endsplit
$$

Thus
$$
\int_{M_T(z_0)}{u^4\over t_0-t}\,do\le
C\Big(E\big(u;D(T;z_0)\big)+
E^{2/3}\big(u;D(T;z_0)\big)\Big).\tag 3.6$$ 
With the special choice $T=0$,
(3.4--6) now lead immediately to Rauch's existence
result: 
\thm{Theorem 3.3}
There exists a constant
$\EPS_0>0$
such that
{\rm(1.8), (1.9)} for any
$u_0\in C^3(\bbbr^3)$,
$u_1\in C^2(\bbbr^3)$ with energy 
$$
E_0=\int_{\bbbr^3}\left({|u_1|^2+|\nabla
u_0|^2\over 2}+{|u_0|^6\over 6}\right)dx<\EPS_0$$ admits 
a global $C^2$-solution.
\ethm
\rem{Remark \rm3.4}
Estimates (3.2)--(3.6) also give the following local 
version of
Rauch's theorem. Let
$\bar z=(\bar x,\bar t)$, $\bar t>0$:
\endrem

\it If $u\in C^2\big(\bbbr^3\times[0,\bar t[\big)$ is
a solution to {\rm(1.8), (1.9)}, and if
$$E\big(u;D(T;\bar z)\big)<\EPS_0\tag3.7$$
for some $T<\bar t$, then $u$ \RM(and its first and second 
derivatives\RM) can be uniformly
a priori bounded on $K_0(\bar z)\smallsetminus\{\bar z\}$ 
in terms of $T,\bar z$,
and the initial data. 

\rm In fact, in this case $u$ can be extended as a solution
of~(1.8) to a full neighborhood
of $\bar z$. Indeed, since $u\in 
C^2\big(\bbbr^3\times[0,\bar t[\big)$,
condition~(3.7) will be satisfied for all points $\tilde 
z=(\tilde x,\bar t)$ with $\tilde
x$ close to $\bar x$.

Finally, observe that if $u$ (and hence its first and 
second derivatives) are uniformly
bounded on $K_0(\bar z)\smallsetminus\{\bar z\}$, 
condition~(3.7) is automatically
satisfied. Thus, in order to extend $u$ as a solution 
of~(1.8) to a neighborhood of
$\bar z$ it suffices to establish that
$$\limsup_{z_0\in K(\bar z)\smallsetminus\{\bar z\}\atop 
z_0\to\bar
z}\big|u(z_0)\big|<\infty.$$
By steps (3.2)--(3.6) of the proof of Rauch's theorem 
then, in fact, it suffices to show
that for some $T<\bar t$ there holds
$$\limsup_{z_0\in K(\bar z)\smallsetminus\{\bar z\}\atop 
z_0\to\bar
z}\int_{M_T(z_0)}{u^4\over t_0-t}do\le 2\pi.\tag3.8$$
This will be important for our next topic.
\heading 4. Large data\endheading
We now show how the smallness assumption in Rauch's
Theorem 3.3 can be removed.
Again remark that it suffices to consider data $u_0,u_1$ 
with compact support.

Suppose by contradiction that $u$ does not extend 
globally. Then $u$ becomes
unbounded in finite time $T$. Since the support of $u$ in 
$\bbbr^3\times[0,T]$ is
relatively compact there exists a ``first" singular point 
$\bar z=(\bar x,\bar t)$,
$0<\bar t\le T$, such that
$$\big|u(x,t)\big|\lra\infty$$
for some sequence $x\to\bar x$,
$t\nearrow\bar t$, and $\bar t$ is minimal with this
property. 

By Remark 3.4, in order to achieve a contradiction it
suffices to establish condition~(3.8). 

Since $t=0$ in the following no longer plays a 
distinguished role, we
may shift coordinates so that $\bar z=(0,0)$
and henceforth assume that
$u$
is a
$C^2$-solution of (1.8) on
$\bbbr^3\times[t_1,0[$
for some
$t_1<0$.
As customary, the Landau symbol ``$o(1)$ as $r\to 0$" 
will denote error terms
depending on a parameter $r$ that tend to 0 as $r\to 0$. 

Observe that (1.8) and $E$ are invariant under scaling
$$
R\mapsto u_R(x,t)=R^{1/2} u(Rx,Rt),
$$
for any $R>0$.

Following Struwe [18, Lemma 2.3] we use the testing 
function
$t\th u_t+x\cdot\nabla u+u$
to derive the following identity
$$0=\bigl(u_{tt}-\2 u+u^5\bigr)\bigl(t\th u_t+
x\cdot\nabla u+u\bigr)
={d\over dt}\bigl(t\th Q_0+u_tu\bigr)-\operatorname{div}
(t\th P_0)+R_0,\tag4.1$$
where
$$\eqalign{P_0
&={x\over t}\left(\tfrac12
|u_t|^2-\tfrac12|\nabla u|^2-\tfrac16|u|^6\right)+
\left(u_t+{x\over t}\cdot\nabla u+{u\over t}\right)\nabla 
u,\cr
Q_0&=\tfrac12|u_t|^2+\tfrac12
|\nabla u|^2+\tfrac16|u|^6+\left({x\over
t}\cdot\nabla u\right)u_t\ge 0\quad \roman{in}\e 
K_{t_1},\cr
R_0&=\tfrac13|u|^6\ge 0.\cr}$$
Note that the multiplier $tu_t+x\cdot\nabla u+u$ is 
related to
the generator of the scaled family $u_R$.
As in Grillakis [6, (2.2)], we may rewrite (4.1) in the 
form
$$0 =t\Biggr\{{d\over dt}\left(Q_0+{u_tu\over t}+{u^2\over
2t^2}\right)-\operatorname{div}
 P_0\Biggr\}+\Biggl\{Q_0+{u^2\over
t^2}+R_0\Biggr\}.\tag4.2$$ 
If we integrate (4.1) over a truncated cone
$K^s_t$,
integrals involving
$u_tu$
will appear. Using the function
$\bigl(t^2+|x|^2\bigr)u_t+2t\th x\cdot\nabla u+2tu$ 
as a further multiplier, Grillakis
succeeds in showing that 
$$
{1\over|t|}\int_{D(t)}u_tu\, dx\le {o}(1)\to
{0}\qquad(t\nearrow 0).$$ 
With little more extra work this additional
multiplier can be avoided. 

As a first step, we obtain
\thm{Lemma 4.1
\rm[18, Lemma 3.2]}
There exists a sequence of numbers
$t_\ell\nearrow 0$
such that
$$
{1\over|t_l|}\int_{D(t_l)}u_tu\th dx\le {o}(1),$$
where
${o}(1)\to 0$
as
$l\to\infty$.
\ethm

For completeness we present the proof.
\demo{Proof}
Consider
$u_l(x,t)=2^{-l/2}u(2^{-l}x,2^{-l}t)$, $l\in\Bbb N$,
satisfying (1.8) with
$$
E\bigl(u_l;D(t)\bigr)=E\bigl(u;D(2^{-l}t)\bigr)\le E_0;$$
moreover, (3.1) translates into the condition
$\bigl( d(u):=d_0(u)\bigr)$
$$\int_{M_{t_1}}d(u_l)\, do\to 0\tag4.3$$
as
$l\to\infty$.
First, suppose that
$$\int_{D(t_1)}u^2_l\, dx\to 0\qquad(l\to\infty).$$
Then let
$t_l=2^{-\ell}t_1$
and estimate
$$
\split
{1\over|t_l|}\int_{D(t_l)}u_tu\,
dx&\le\left(\int_{D(t_l)}|u_t|^2\,dx\right)^{1/2}
\left({1\over|t_l|^2}\int_{D(ta_l)}u^2\,dx\right)^{1/2}\\
&\le\Bigl(2E\bigl(u;D(t_l)\bigr)\Bigr)^{1/2}
\left({1\over|t_1|^2}\int_{D(t_1)}u^2_l\,
dx\right)^{1/2}\to 0\quad(l\to\infty)
\endsplit
\tag4.4
$$
to achieve the claim.

Otherwise, there exist
$C_1>0$
and a sequence
$\Lambda$
of numbers
$l\to\infty$
such that
$$\liminf_{l\to\infty,\th l\in\Lambda}\int
_{D(t_1)}u^2_l\, dx\ge C_1.$$
For any
$s\in[t_1,0[$,
by H\"older's inequality
$$\int_{D(s)}u^2_l\,
dx\le\left({4\pi\over3}|s|^3\right)^{2/3}\left(\int
_{D(s)}|u_l|^6\,dx\right)^{1/3}\le
C\th E_0^{1/3}s^2.\tag4.5$$ Choose $s=s_1<0$
such that the latter is
$\le C_1$.
Then by (4.3) we have
$$\eqalign{2\int_{K_{t_1}^{s_1}}(u_l)_tu_l\, dz&
=\int_{K_{t_1}^{s_1}}{d\over
dt}|u_l|^2\,dx\cr
&=\int_{D(s_1)}|u_l|^2\,dx-\int_{D(t_1)}|u_l|^2\,dx+
{1\over\sqrt
2}\int_{M_{t_1}^{s_1}}|u_l|^2\,do\cr
&\le {o}(1)\to 0\qquad(l\to\infty,\e l\in\Lambda).\cr}$$
We conclude that for suitable numbers
$s_l
\in[t_1,s_1]$,
$t_l=2^{-\ell} s_l$, $l\in\Lambda$,
we have
$${2\over|t_l|}\int
_{D(t_l)}u_tu\, 
dx={2\over|s_\ell|}\int_{D(s_l)}(u_l)_tu_l\,
dx\le {o}(1)\to 0\qquad(l\to\infty,\ l\in\Lambda).$$
Relabelling, we obtain a sequence
$(t_l)_{l\in\Bbb N}$,
as desired.\qed\enddemo
\thm{Lemma 4.2 \rm[18, Lemma 2.2]}
For any
$l\in\Bbb N$
there holds
\rm
$$
\split
&{1\over 3|t_l|}\int_{K_{t_l}}|u|^6\,dz+\int
 _{D(t_l)}\left({1\over2}|u_t|^2+{1\over 2}|\nabla u|^2
 +u_t\left({x\over t}\cdot\nabla u\right)+{1\over
6}|u|^6\right)\,dx\\
&\qquad\le{o}(1)\to 0\qquad(l\to\infty).
\endsplit
$$
\ethm

\n
Again we give the proof for completeness.
\demo{Proof}
For
$s\in[t_l,0[$
integrate (4.1) over
$K^s_{t_l}$
to obtain
$$\split
0=&\int_{D(s)}\bigl( s \th Q_0+u_tu\bigr)\,dx+{1\over\sqrt
2}\int_{M^s_{t_l}}\bigl(t\th Q_0+u_tu-x\cdot 
P_0\bigr)\,do\\
&+|t_l|\int_{D(t_l)}Q_0\,dx+
\int_{K^s_{t_l}}R_0\,dx-\int_{D(t_l)}u_tu\, dx.
\endsplit
$$
Now,
$Q_0$
is dominated by the energy density
$e(u)$.
Thus, and using H\"older's inequality as in (4.4), (4.5), 
the first term is of order
$|s|$
and hence vanishes as
$s\nearrow 0$.
Moreover, on
$M_{t_l}$
we have
$$
\split
&\bigl| t\th Q_0+u_tu-x\cdot P_0\bigr|\\
&\qquad=|t|\left| |\nabla
u|^2-\left|{x\over|x|}\cdot\nabla u\right|^2
 +\tfrac13 |u|^6-{u(t\th u_t+x\cdot\nabla u)\over 
t^2}\right|\\
&\qquad\le|t_l|\left(3d_0(u)+{|u|^2\over 
t^2}\right).\endsplit
$$
Hence by (3.1) and Hardy's inequality Lemma 3.1 the 
second term is of order
${o}(1)|t_l|$,
where
${o}(1)\to 0$
as
$l\to\infty$.
Thus, by Lemma 4.1 we have
$$
{1\over|t_l|}\int_{K^s_{t_\ell}}R_0\,dz+
\int_{D(t_l)}Q_0\,dx\le{1\over
|t_l|}\int_{D(t_l)}u_tu\, dx+{o}(1)\le{o}(1)\to 0
\qquad(l\to\infty),
$$
which is the desired conclusion.\qed\enddemo

Now we use (4.1) in its equivalent form (4.2) to derive a 
stronger version of
Lemma 4.1.  
\thm{Lemma 4.3}
There exists a sequence of numbers
$\bar t_l\nearrow 0$
such that the conclusion of Lemma {\rm4.1} holds for
$(\bar t_l)$
while in addition we have
$$2\le \bar t_l{\di /}\bar t_{l+1}\le 4$$
for all
$l\in\Bbb N$.
\ethm
\demo{Proof}
First observe that by H\"older's inequality and Lemma 4.2 
for any
$m\in\Bbb N$
we have
$$\int_{D(t_m)}{|u|^2\over |t|^2}\,dx\le C\left(\int
_{D(t_m)}|u|^6\,dx\right)^{1/3}\le C
\left(\int_{D(t_m)}Q_0 \,dx\right)^{1/3}\longrightarrow 
0\quad(m\to\infty),$$
where
$(t_m)$
is determined in Lemma 4.1.
Divide (4.2) by
$t$
and integrate over the cone
$K^{t_m}_{t_l}$
for 
$m\ge l$
to obtain
$$
\split
&\int_{D(t_l)}\left(Q_0+{u_tu\over t}+{u^2\over
2t^2}\right)dx+\int_{K^{t_m}_{t_l}}
\left({Q_0\over|t|}+{|u|^2\over |t|^3}+{R_0\over
|t|}\right)dz\\
&\qquad=\int_{D(t_m)}\left(Q_0+{u_tu\over t}+{u^2\over
2t^2}\right)dx
+\int_{M^{t_m}_{t_l}}\left(Q_0+{u_tu\over t}+{u^2\over 
2t^2}-{x\over
t}\cdot P_0\right)do.\endsplit
$$
By the preceding remark the first term on the right 
vanishes as we let
$m\to\infty$,
while by (3.1) the integral over
$M^{t_m}_{t_\ell}$
becomes arbitrarily small as
$m\ge l\to\infty$.
Finally, by Lemma 4.1, we have
$$\int_{D(t_l)}{u_tu\over t}\,dx=-{1\over|t_l|}\int 
u_tu\, dx\ge {o}(1)\to
0\qquad (l\to\infty).$$
All remaining terms being nonnegative, we thus obtain the 
estimate
$$
\int
_{K_{t_l}}{|u|^2\over |t|^3}\,dz=\int^0_{t_l}
\left({1\over |t|}\int_{D(t)}{|u|^2\over
|t|^2}\,dx\right)\,dt\le {o}(1)\to 0\qquad(l\to\infty).$$
Hence for any
$\bar t\in[t_l/2,0[$
there also holds
$$
{o}(1)\ge{1\over \bar t}\int^{\bar t}_{2\bar
t}
\left(\int
_{D(t)}{|u|^2\over|t|^2}\,dx\right)dt\ge \inf_{2\bar t\le 
t\le\bar
t}\int_{D(t)}{|u|^2\over |t|^2}\,dt,$$ where
${o}(1)\to 0$
if
$l\to\infty$.

To construct the sequence
$(\bar t_l)$,
now choose
$\bar t_1=t_1$
and proceed by induction. Suppose
$\bar t_l$, $l=1,...,L$,
have been defined already.
Let
$\bar t_{L+1}\in\left[{\bar t_L\over 2},{\bar t_L\over 
4}\right[$
be chosen such that
$$
\int_{D\bigl(\bar t_{L+1}\bigr)}
{|u|^2\over |t|^2}\,dx\le2\inf_{{\bar t_L\over 2}\le t\le
{\bar t_L\over 4}}\int_{D(t)}{|u|^2\over|t|^2}\,dx.$$
Clearly,
this procedure yields a sequence
$(\bar t_l)$
such that
$2\le \bar t_l /\bar t_{l+1}\le 4$
for all
$l$ and
$$\int_{D(\bar t_l)}{|u|^2\over|t|^2}\,dx\longrightarrow 
0\qquad
(l\to\infty).$$

By (4.4) then
$$
{1\over |\bar t_l|}\int_{D(\bar t_\ell)}u_tu\, 
dx\longrightarrow
0\qquad(l\to\infty),$$ concluding the proof.\qed\enddemo

\n
To simplify notation, we will assume that
$t_l=\bar t_l$
for all
$l$,
initially.
\subheading{The radial case}
At this point we can indicate how the decay estimate 
Lemma 4.2 and Lemma 4.3
imply regularity of solutions in the radial case. First 
observe that for radially
symmetric data $u_0(x)=u_0\bigl(|x|\bigr)$, etc., the 
unique local 
$C^2$-solution
$u$
to (1.8), (1.9) again is radially symmetric, that is,
$u(x,t)=u\bigl(|x|,t\bigr)$. 

Note that this implies that blow-up can only occur on
the line $x=0$. Indeed, if $u$
is regular on
$K_0(\bar z)\backslash\{\bar z\}$
and blows up at
$\bar z=(\bar z, \bar t)$,
the same will be true for any point
$z=(x,\bar t)$
with
$|x|=|\bar x|=\bar r$. Now, if
$\bar x\ne 0$,
given any
$K\in\Bbb N$
we can choose points
$x_k\in\bbbr^3$ with $|x_k|=\bar r$,
$1\le k\le K$, and $T\in[0,\bar t[$
such that
$$
D(T;z_k)\cap D(T;z_l)=\emptyset$$
for all
$k\ne l$, where $z_k=(x_k,T),\ k=1,...,K$.
Moreover, by Remark 3.4
$$
E\bigl(u;D(T;z_k)\bigr)\ge \EPS_0>0$$
for all
$k$, while by Lemma 3.2
$$\eqalign{K\EPS_0&\le\sum^K_{k=1}E\bigl(u;D(T;z_k)\bigr)
=E\left(u;\bigcup^K_{k=1}D(T;z_k)\right)\le 
E\Bigl(u;D\bigl(T;(0, \bar r+\bar
t)\bigr)\Bigr)\cr
&\le E\Bigl(u;D\bigl(0;(0,\bar r+\bar t)\bigr)\Bigr)\le 
E_0,\cr}$$ 
independently of
$K$. Thus, blow-up may first occur on the line
$x=0$,
only. (See Figure~1.)

Let blow-up occur at
$(0,\bar t)$. Shifting time by $\bar t$
then we may assume that
$u(x,t)=u\bigl(|x|,t\bigr)$
is regular on
$\bbbr^3\times[t_1,0[$
and blows up at the origin. 
As a second step we estimate the speed of blow-up. Again 
observe that (1.8) is
invariant under scaling $$u\longmapsto 
u_R(x,t)=R^{1/2}u(R\th x,R\th t).$$ 
This suggests that
$u(z)\sim|z|^{-1/2}$. In fact, the following result holds.

\fighere{17pc}\caption{\smc Figure 1. \rm The energy
at the basis of each cone is $\ge\varepsilon_0$.}
\eject

\topspace{18pc}\caption{{\smc Figure 2. \rm Overlap of 
the cones $K
_{t_{l-L}}(z^j_k)$.}}

\thm{Lemma 4.4 \rm[18; Lemma 3.3]}
$$4\EPS_1:=\limsup_{z=(x,t)\to 0\atop z\in
K_{t_1}}\Bigl(\bigl|u(z)\bigr|\cdot|z|^{1/2}\Bigr)>0.$$
\ethm

The proof of Lemma 4.4 is rather involved and will not be 
presented here.

We can now conclude the regularity proof. Let
$\bar z_k=(\bar x_k,\bar s_k)\in K_{t_1}$
satisfy
$$\bigl|u(\bar z_k)\bigr|=\sup_{z\in K_{t_1}(\bar 
z_k)}\bigl| u(z)\bigr|\ge 2\EPS_1|\bar
z_k|^{-{1\over 2}}.$$
Choose a sequence
$l=l(k)\to \infty\e(k\to\infty)$
such that
$$t_{l+1}\ge\bar s_k\ge t_l,$$
with
$(t_l)$
as in Lemma 4.3.
By (3.2)--(3.3) we may fix $L\in\Bbb N$
independent of
$k$
such that for large
$k$
there holds
$$\align
\bigl|u(\bar 
z_k)\bigr|\left(1-{1\over4\pi}\int_{M_{t_{l-L}}
(\bar z_k)}{u^4\over \bar s_k-t}\,do\right)
&\le C+C\bigl(|\bar s_k-t_{l-L}|^{-1/ 2}-|\bar
s_k-t_1|^{-1/2}\bigr)E_0^{5/6}\\
&\le C+C\th 2^{-L/2}|\bar 
s_k|^{-{1/2}}E_0^{5/6}\le\EPS_1|\bar z_k|
^{-{1/2}},\endalign$$
where
$E_0=E\bigl(u;D(t_1,0)\bigr)$
is the initial energy.

Thus, by choice of
$\bar z_k$
and
(3.6) we conclude that
$$E\bigl(u;D(t_{l-L};\bar z_k)\bigr)\ge\EPS_0\quad
\text{for all }k.$$
Given
$J\in\Bbb N$,
for each
$k\in\Bbb N$ choose 
$J$
points
$x^j_k$, $j=1,...,J$,
equi-distributed on
the sphere
$|x_k^j|=|\bar x_k|$. Let
$z^j_k=(x_k^j,\bar s_k)$. Note that there exists
$\dt=\dt(J,L)>0$
such that
$(x,t_{l-L})\in D(t_{l-L};z^i_k)\cap D(t_{l-L}; z^j_k)$, 
$i\ne j$,
implies that
$|x|\le(1-\dt)t_{l-L}$. (See Figure~2.)

Hence by
Lemma 4.2 we have
$$\eqalign{J\EPS_0&\le\sum^J_{j=1}E\bigl(u;D(t_{l-L}; 
z_k^j)\bigr)\cr
& \le E\bigl(u;D(t_{l-L};0)\bigr)+\sum_{i\ne 
j}E\bigl(u;D(t_{l-L};z^i_k)\cap
D(t_{l-L};z^j_k)\bigr)\cr
&\le E_0+C(J,\delta)\int_{D(t_{l-L})} Q_0\,dx\le E_0+
o(1),\cr}$$
where
$o(1)\to0$
as
$k\to\infty$
for any fixed
$J$.
Choosing
$J$
large, for
sufficiently large
$k\in\Bbb N$
we thus obtain a contradiction. Hence,  $u$ is uniformly 
bounded on
$K_{t_1}$ and the proof in the radially symmetric case is 
complete.
\subheading{General data}
Finally, we present Grillakis' work on the general case. 
The key observation is that
the decay Lemma 4.2 suffices to establish (3.8), 
directly. However, this is not at all
easy to see.

Fix
$z_0=(x_0,t_0)\in K\backslash\{0\}$
arbitrarily. Denote
$y=x-x_0$, $\hat y={y\over|y|}$, $\hat x={x\over|x|}$.
Divide (4.2) by
$t$
and for
$s>t_0$ 
integrate over
$K^s_{t_l}\backslash K(z_0)$
to obtain the relation
$$
\split
0=&\int_{D(s)}\left(Q_0+{u_tu\over t}+{u^2\over
2t^2}\right)\,dx
-\int_{D(t_l)\backslash D(t_l;z_0)}\left(Q_0+{u_tu\over
t}+{u^2\over 2t^2}\right)dx\\
&+{1\over\sqrt 2}\int_{M^s_{t_l}}\left(Q_0+{u_tu\over t}+
{u^2\over
2t^2}-\hat x\cdot P\right)do\\
&-{1\over\sqrt 2}\int_{M_{t_l}(z_0)}\left(Q_0+{u_tu\over 
t}+{u^2\over
2t^2}-\hat y\cdot P\right)do\\
&+\int_{K^s_{t_l}\backslash K(z_0)}\left({R_0+Q_0+
{u^2\over 2t^2}\over
t}\right)dz=I+\dots+V.\endsplit
$$
(See Figure 3.)

\fighere{16.5pc}\caption{\smc Figure 3}
\eject
By H\"older's inequality (4.4), (4.5) and Lemma 4.2 the 
first term
I $\to0$
if we choose
$s=t_k$
with
$k\to\infty$.
Similarly,
II $\to 0$
if
$l\to\infty$.
By (3.1) also III $\to 0$
as
$l\to\infty$. Finally
V
$\le 0$.
Thus we obtain the estimate for any
$z_0\in K\backslash\{0\}$,
$$
\int_{M_{t_l}(z_0)}\left(Q_0+{u_tu\over t}+{u^2\over 
2t^2}-\hat y\cdot
P\right)do\le {o}(1)\to 0\qquad(l\to\infty),\tag4.6$$
with error term ${o}(1)$
independent of $z_0$.

By a beautiful geometric reasoning, Grillakis [6] now 
proceeds to bound (3.8) in
terms of (4.6) of (4.6). Let
$r=|x|$;
then we may rewrite
$$\eqalign{A:&=Q_0+{u_tu\over t}+{u^2\over 2t^2}-\hat 
y\cdot P\cr
&=\left(1-\hat x\cdot\hat y{r\over 
t}\right)\tfrac12|u_t|^2+
\left(1+\hat x\cdot\hat y{r\over 
t}\right)\left(\tfrac12|\nabla u|^2+
\tfrac16 |u|^6\right)\cr
&\quad+{1\over t}\bigl(u_t-\hat y\cdot\nabla u\bigr)u+
{r\over t}u_t\hat x\cdot\nabla
u-u_t\hat y\cdot\nabla u
-{r\over t}\bigl(\hat x\cdot\nabla u\bigr)\bigl(\hat 
y\cdot\nabla
u\bigr)+{u^2\over 2t^2}.\cr}$$
Introducing
$u_\sigma=\hat y\cdot\nabla u,\th\a=\hat x-\hat y(\hat 
y\cdot\hat x)$,
$|\a|u_\a=\a\cdot\nabla u,\th\1 u=\nabla u-\hat y 
u_\sigma$,
this expression becomes
$$
\eqalign{=&\left(1-\hat x\cdot\hat y{r\over 
t}\right)\tfrac1
2(u_t-u_\sigma)^2+\left(1+\hat x\cdot\hat y{r\over 
t}\right)\left(\tfrac
12|\1
u|^2+\tfrac16|u|^6\right)\cr
&+{r\over t}|\a|(u_t-u_\sigma)u_\a+{u\over 
t}(u_t-u_\sigma)+{u^2\over 2t^2}.\cr}$$
Now let
$\hat x\cdot\hat y=$
cos
$\delta$,
$|\a|=$
sin
$\de$
and for
brevity denote
${1\over\sqrt 2}(u_t-u_\sigma)=\ur$. (See Figure~4 on p\. 
78.)
Then the above

\topspace{23.5pc}\caption{\smc Figure 4} 

$$
\split
A=&\left(1-{r\over 
t}\th\operatorname{cos}\th\delta\right)|\ur|^2+
\left(1+{r\over 
t}\th\operatorname{cos}\th\delta\right)\left(\tfrac12
|\1 u|^2+\tfrac16|u|^6\right)\\
&+\sqrt 2{r\over t}|\sin\th\delta|\ur u_\a+\sqrt 2{\ur 
u\over t}+{u^2\over
2t^2}\\
=&:A_0+\sqrt 2{\ur u\over t}+{u^2\over 2t^2}.\endsplit
\tag4.7$$

\n
Note that if we estimate
$|u_\a|\le|\1 u|$,
in particular, we have
$$\split
A_0\ge&\left(1-{r\over t}\th\cos \th\delta\right)|\ur|^2+
\left(1+{r\over 
t}\th\cos\th\delta\right)\left(\tfrac12|u_\a|^2+\tfrac16
|u|^6\right)
+\sqrt 2{r\over t}|\sin\th\de|\ur u_\a\\
\ge&\left(1+{r\over t}\right)\!\left(|u_\rho|^2+\frac12
u_\a|^2+\tfrac16|u|^6\right)\!-\!{r\over 2t}\Bigl(\sqrt 
2\sqrt{1+\cos
\th\de}\ur-\sqrt{1-\cos\th\de}u_\a\Bigr)^2\\\ge& 0\endsplit
\tag4.8$$
on
$M_{t_l}(z_0)$.

Now for any 
$\EPS >0$
there is a constant
$C=C(\EPS)$
such that for any
$z_0\in K$
and any
$z\in M^{Ct_0}(z_0)$
we may estimate\rm
$$
-\sqrt 2{r\over t}\th|\sin\th\de|\le\EPS,\qquad
-{r\over t}\th\cos\th
\de\ge\tfrac12.$$
In fact, for
$z=(x,t)\in M^{Ct_0}(z_0)$
we have
$$\bigl| |x|-|y|\bigr|\le |y-x|=|x_0|\le 
|t_0|\le{|t-t_0|\over C-1}={|y|\over C-1}.$$

Hence
$$\hat x\cdot\hat y=\th\cos\th\de={x\over 
|x|}\cdot{y\over |y|}\ge
1-\left|{y\over |y|}-{x\over |x|}\right|\ge 
1-2{|x_0|\over |y|}\ge 1-{2\over C-1}$$ 
while
$$1\ge-{r\over t}={|y|\over|t-t_0|}\cdot{|t-t_0|\over 
|t|}\cdot{|x|\over
|y|}\ge\left(1-{1\over C}\right)\left(1-{1\over 
C-1}\right).$$
This yields the following estimate.
\thm{Lemma 4.5}
For any
$\EPS>0$,
any
$z_0\in K$,
letting
$C=C(\EPS)$
be determined as above, if
$t_k\le C\th t_0$
we have
$$\int_{M^{t_k}_{t_l}(z_0)}A \, do\ge{1\over
2}\int_{M^{t_k}_{t_l}(z_0)}|\ur|^2\,do-\EPS E_0.$$
\ethm
\demo{Proof}
Estimate
$$\left|{\sqrt 2\ur u\over t}\right|\le|\ur|^2+{u^2\over 
2t^2}.$$
Hence by (4.7) and our choice of $C(\EPS)$, for
$z\in M^{C\th t_0}_{t_l}(z_0)$
we have
$$A\ge\tfrac12
|\ur|^2-\EPS|\ur u_\a|\ge\tfrac12
|\ur|^2-\tfrac{1}{\sqrt{2}}\in d_{z_0} (u),$$
which in view of Lemma 3.2 proves the claim.\qed\enddemo

Observe that
$\ur$
may be interpreted as a tangent derivative along $M(z_0)$.
In fact, let
$\Phi$
be the map
$$\Phi:y\mapsto(x_0+y,\th t_0-|y|)\tag4.9$$
and let
$$v(y)=u\bigl(\Phi(y)\bigr),\tag4.10$$
whenever the latter is defined. Then the radial derivative
$v_s$
of
$v$
is given by
$$v_s={y\cdot\nabla v\over |y|}=u_\sigma-u_t=-\sqrt 2 
\ur.\tag4.11$$
\thm{Lemma 4.6 \rm(Grillakis [6, (2.23)])}
For any
$z_0\in K$
and any
$C\ge 0$
there holds
$$
\int_{M_{(1+C)t_0}(z_0)}{\ur u\over t}\,do\ge\bigl(1+\ln
(1+C)\bigr)\cdot {o}(1)$$
where
${o}(1)\to 0$
if
$(1+C)t_0\ge t_l$
and
$l\to\infty$.
\ethm
\demo{Proof}
Introducing
$y$
as new variable via (4.9), (4.10), we have
$$
\int_{M_{(1+C)t_0}(z_0)}{\ur u\over t}\,do=\int
_{B_{C|t_0|}}v_s{v\over
|y|-t_0}dy
=\int_{S_1}\left(\int
^{C|t_0|}_0v_s{v\over s-t_0}s^2\,ds\right)\,do.
$$
Upon integrating by parts, this gives
$$\eqalign{\int_{S_1}&\left(\int^{C|t_0|}_0{\partial_s(v 
^{2}/2)\over
s-t_0}s^2\,ds\right)\,do\cr
&=\int_{S_1}\int^{C|t_0|}_0\left(-{v^2s\over s-t_0}+
{v^2s^2\over
2(s-t_0)^2}\right)ds\, do
+{1\over 2(1+C)|t_0|}\int_{S_{C|t_0|}}v^2\,do\cr
&\ge-\int_{B_{C|t_0|}}{v^2\over 
|y|\bigl(|y|-t_0\bigr)}\,dy=-{1\over \sqrt
2}\int_{M_{(1+C)t_0}(z_0)}{u^2\over t(t - t_0)}\,do(x,t)\cr
&=-\int^{t_0}_{(1+C)t_0}{1\over 
|t|}\left({1\over|t-t_0|}\int_{\partial
D(t;z_0)}u^2\,do(x)\right)\,dt.\cr}$$
Now by Hardy's inequality Lemma 3.1.(iii) we have
$$\eqalign{&\left({1\over |t-t_0|}\int_{\partial D(t;z_0)}
\negthinspace\negthinspace\negthinspace\negthinspace%
\negthinspace\negthinspace\negthinspace 
u^2\,do\right)^2
\le C\int_{\partial D(t;z_0)}
\negthinspace\negthinspace\negthinspace\negthinspace%
\negthinspace\negthinspace\negthinspace 
u^4\,do\cr 
&\qquad\le C\left\{\left(\int_{D(t;z_0)}
\negthinspace\negthinspace\negthinspace\negthinspace%
\negthinspace\negthinspace\negthinspace 
|\nabla
u|^2\,dx\right)^{1/2}+ \left(\int_{
D(t;z_0)}
\negthinspace\negthinspace\negthinspace\negthinspace%
\negthinspace\negthinspace\negthinspace
|u|^6\,dx\right)^{1/6}\right\}
\cdot\left(\int_{
D(t;z_0)}
\negthinspace\negthinspace\negthinspace\negthinspace%
\negthinspace\negthinspace\negthinspace
|u|^6\,dx\right)^{1/2}\cr 
&\qquad\le C(E_0)\left(\int_{
D(t;z_0)}
\negthinspace\negthinspace\negthinspace\negthinspace%
\negthinspace\negthinspace\negthinspace
|u|^6\,dx\right)^{1/2}.\cr}$$
Hence 
$$\int_{M_{(1+C)t_0}(z_0)}{\ur u\over t}\, do\ge - 
C\int^{t_0}_{(1+C)t_0}{1\over
|t|}\left(\int_{D(t)}|u|^6\,dx\right)^{1/4}\,dt,$$
with $C$
depending on
$E_0$
only. By Lemma 4.2 the latter can be controlled as
follows. Let
$k,K\in\Bbb N$
be determined such that
$$t_k\le(1+C)t_0<t_{k+1}\le t_{K}\le t_0<t_{K+1}.$$

Note that by Lemma 4.3
$$1+C\ge{t_{k+1}\over {t_K}}\ge 2^{K-(k+1)},$$
whence
$$K-k\le 1+\th\log_2(1+C).$$

Estimate
$$\int^{t_0}_{(1+C)t_0}{1\over 
|t|}\left(\int_{D(t)}|u|^6\,dx\right)^{1/4}dt\le
\sum^K_{i=k}\int^{t_{i+
1}}_{t_i}{1\over|t|}\left(\int_{D(t)}|u|^6\,dx%
\right)^{1/4}dt.$$
By H\"older's inequality, this is
$$\le C\sum^K_{i=k}{\bigl|
t_i-t_{i+1}\bigr|^{3/4}\over\bigl|t_{i+1}\bigr|}
\left(\int_{K^{t_{i+1}}_{t_i}}|u|^6\,dz\right)^{1/4}$$
and by Lemma 4.3
$$\le C\sum^K_{i=k}\left({1\over 
|t_i|}\int_{K_{t_i}}|u|^6\,dz\right)^{1/4}.$$
Finally, use Lemma 4.2 to see that this is
$$\le(K-k+1){o}(1)\leqq \bigl(1+\ln(1+C)\bigr){o}(1)$$
where
${o}(1)\to$
if
$(1+C)t_0\ge t_l$
and
$l\to\infty$.\qed\enddemo

Combining Lemmas 4.5 and 4.6 it follows that for any
$\EPS>0$,
if we choose
$t_k\le C(\EPS)t_0<t_{k+1}$,
we can estimate
$$\eqalign{{o}(1)&\ge\int_{M_{t_l}(z_0)}A\th do\cr
&\ge\tfrac12\int_{M^{t_k}_{t_l}(z_0)}|\ur|^2\,do-\EPS E_0
+\int_{M_{t_k}(z_0)}A_0\,do-{o}(1)\Bigl(1+\ln\bigl(1+
C(\EPS)\bigr)\Bigr),\cr}
\tag4.12$$
where ${o}(1)\to0$
as
$l\to\infty$.
To estimate
$A_0$
on
$M_{t_k}(z_0)$
now introduce the new angle
$\de_0$,
where
$|x_0|=r_0,\e\hat x_0={x_0\over r_0},\e\hat x_0\cdot \hat 
y=$
cos
$\de_0$.
(See Figure~5.)
Again
$y=x-x_0$,
and
$|y|=\sigma=t_0-t$.
With this notation

\topspace{22.5pc}\caption{\smc Figure 5}

$$r\th\cos\th\dt=\hat x\cdot\hat y r=x\cdot\hat 
y=y\cdot\hat y+x_0\cdot\hat
y=\sigma+r_0\th\cos\th\de_0,$$
$$
|\sin\th\dt|=|\a|
=\left|{x-(x\cdot\hat y)\hat y\over r}\right|=
\left|{x_0-(x_0\cdot\hat y)\hat y\over r}
\right|={r_0\over r}|\sin\th\de_0|.$$

\n
Hence, by (4.7),
$$\eqalign{A_0&=\left(1-{\sigma\over t}-{r_0\over t}\th\cos
\th\de_0\right)|\ur|^2
+\left(1+{\sigma\over t}+{r_0\over t}\th\cos
\th\de_0\right)\left(\tfrac12|\1 u|^2+
\tfrac16|u|^6\right)\cr
&\quad+\sqrt 2\e {r_0\over t}|\sin\th\de_0|\ur 
u_\a.\cr}\tag4.13$$
Estimating
$|\1 u|\ge|u_\a|$
as before, this is
$$\eqalign{&\ge\left(2-{t_0-r_0\over t}\right)|\ur|^2
-{r_0\over 2t}\left(\sqrt 2\sqrt{1+
\th\cos\th\de_0}\ur-\sqrt{1-\th\cos
\th\de_0}u_\a\right)^2\cr
&\quad+{t_0\over t}\left(1+{r_0\over t_0}\right)\tfrac12
|u_\a|^2+{t_0\over
t}\left(1+{r_0\over t_0}\th\cos\th\de_0\right)\tfrac16
|u|^6.\cr}\tag4.14$$
Note that all the latter terms are nonnegative for
$z\in M(z_0),\th z_0\in K$.
By (4.14),
and since
$r_0\le|t_0|$,
for
$t\le 2t_0$
we have
$A_0\ge|\ur|^2$.
Moreover, given
$0<\EPS<1,\th
z_0\in K$,
let
$t_m\le 2t_0<t_{m+1}$
and set
$$\eqalign{\Gamma&=\Gamma(\EPS;\th z_0)\th=\left\{z\in 
M_{t_m}(z_0);\e
|\de_0|\le\EPS^{1/4}\right\}\cr
\2&=\2(\EPS;\th z_0)=M_{t_m}(z_0)\backslash\Gamma.\cr}$$

Note that by (4.13) on
$\Gamma$
we can estimate
$$\eqalign{A_0&\ge|\ur|^2-\sqrt 2\EPS^{1/4}|\ur u_\a|\cr
&\ge|\ur|^2-\sqrt 2\EPS^{1/4}d_{z_0}(u),\cr}$$
while by (4.14) on
$\2$
we have
$$\eqalign{A_0&\ge{t_0\over t}\left(1+{r_0\over t_0}\th\cos
\th\de_0\right)\tfrac16|u|^6
\ge\tfrac18
\Biggl(1-\left(1-{\EPS^{1/2}\over 2}+
\EPS\right)\Biggr)\tfrac16
|u|^6\cr
&\ge{\EPS^{1/2}\over 96}|u|^6-\EPS d_{z_0}(u).\cr}$$
Combining with (4.12) and Lemma 4.5, thus we obtain
$$\split
\int
_\Gamma|\ur|^2do&\le\int_{M_{t_k}(z_0)}A_0do+\sqrt 
2\EPS^{1/4}E_0\\
&\le\left(\EPS+\sqrt 2\EPS^{1/4}\right) E_0+o(1)\Bigl(1+
\ln \bigl(
1+C(\EPS)\bigr)\Bigr),\endsplit\tag4.15$$ 
$$
{\EPS^{1/2}\over
96}\int_\2|u|^6\,do\le\int_{M_{t_k}(z_0)}A_0do+\EPS E_0
\le 2\EPS E_0+{o}(1)\Bigl(1+\ln \bigl(
1+C(\EPS)\bigr)\Bigr),\tag4.16$$ 
$$
\int_{M^{t_m}_{t_l}(z_0)}\!|\ur|^2\,do\le\int
_{M^{t_m}_{t_k}(z_0)}\!A_0\,do+\!
\int_{M^{t_k}_{t_l}(z_0)}|\ur|^2\,do
\le 2\varepsilon E_0+{o}(1)\ln \bigl(
1+C(\EPS)\bigr)\!,\tag4.17$$ 
where ${o}(1)\to0$ as $l\to\infty$. (We may assume 
$t_l\le t_k\le t_m$.)
\demo{Proof of Theorem \rm1.1}
Given
$\EPS>0$,
we split the integral in (3.8) and use H\"older's 
inequality as follows
$$
\split
\il_{M_{t_l}(z_0)}{|u|^4\over |t-t_0|}\,do&\le\il_\Gamma 
+\dots+\il_\2 +\dots+
\il_{M^{t_m}_{t_l}(z_0)}+\dotsb\\
&\le\left(\il_\Gamma{|u|^2\over
|t-t_0|^2}do\right)^{1/2}\left(\il_\Gamma|u|^6\,do%
\right)^{1/2}\\
&\quad+\left(\il_\2{|u|^2\over
|t-t_0|^2}\,do\right)^{1/2}\left(\il_\2|u|^6\,do%
\right)^{1/2}\\
&\quad+\left(\il_{M^{t_m}_{t_l}(z_0)}{|u|^2\over
|t-t_0|^2}\,do\right)^{1/2}
\left(\il_{M^{t_m}_{t_l}(z_0)}|u|^6\,do\right)^{1/2}.%
\endsplit
$$
By Lemma 3.1.(ii) and Lemma 3.2 this can be further 
estimated
$$
\eqalign{&\le\sqrt{6E_0}\left(\il_\Gamma {|u|^2\over
|t-t_0|^2}\,do\right)^{1/2}-{C(E_0)}\left(\il_\2|u|^6\,do%
\right)^{1/2}\cr
&\quad+\sqrt{6E_0}
\left(\il_{M^{t_m}_{t_l}(z_0)}{|u|^2\over|t-t_0|^2}\,do%
\right)^{1/2}.\cr}$$
By Lemma 3.1.(i)
and (4.15)
$$\eqalign{\il_\Gamma{|u|^2\over |t-t_0|^2}\,do&\le
4\il_\Gamma|\ur|^2\,do+2|t_m-t_0|^{-1}\il_{\partial 
D(t_m;\th z_0)}|u|^2\,do\cr
&\le 4 \left(\EPS+\sqrt
2\EPS^{1/4}\right)E_0+{o}(1)\Bigl(1+\ln\bigl(1+
C(\EPS)\bigr)\Bigr)\cr
&\quad+C\left(
\il_{\partial D(t_m;\th z_0)}|u|^4\,do\right)^{1/2}.\cr}$$
By Lemma 3.1.(iii) and Lemma 4.2 the latter
$$\align
\il_{\partial D(t_m,z_0)}\negthinspace
\negthinspace 
\negthinspace\negthinspace\negthinspace\negthinspace 
 |u|^4\,do
&\le C\!\left(\left(\il_{D(t_m)}\negthinspace
\negthinspace\negthinspace\negthinspace|\nabla
u|^2\,dx\right)^{1/2}\!+\!\left(\il_{D(t_m)}
\!|u|^6\,dx\right)^{1/6}\right)\!\left(\il_{D(t_m)}%
\negthinspace
\negthinspace\negthinspace\negthinspace 
|u|^6\,dx\right)^{1/2}\\
&\le C(E_0){o}(1),\endalign
$$
where
${o}(1)\to0$
as
$m\ge l$
tend to infinity.
Similarly, by Lemma 3.1.(i), (iii), Lemma 4.2 and (4.17)
$$
\align
\il_{M^{t_m}_{t_l}(z_0)}{|u|^2\over |t-t_0|^2}\,do&\le 4
\il_{M^{t_m}_{t_l}(z_0)}|\ur|^2\,do +
2|t_l-t_0|^{-1}\il_{\partial
D(t_l;z_0)}|u|^2\,do\\
&\le \EPS
E_0+{o}(1)\Bigl(C(E_0)+\ln\bigl(1+
C(\EPS)\bigr)\Bigr).\endalign$$ 

Finally, by (4.16)
$$\il_\2|u|^6\,do\le
192\EPS^{1/2}E_0+{o}(1)\EPS^{-1/2}\Bigl(1+\ln\bigl(1+
C(\EPS)\bigr)\Bigr).$$
Hence, if we first choose
$\EPS>0$
sufficiently small and then choose
$l\in\Bbb N$
sufficiently large, the integral
$$\il_{M_{t_l}(z_0)}{u^4\over |t-t_0|}\,do$$
can be made as small as we please, uniformly in $z_0\in 
K_{t_1}$.\qed\enddemo
\rem{Remark}
Since all error estimates are based on the qualitative 
statement (3.1), no a priori
bounds for the solution 
$u$
on a cone
$K$, depending only on
$u_0,u_1$,
and
$K$,
are obtained.
\endrem
\heading 5. A remark on the super-critical case\endheading
We add an observation on the
super-critical case. Consider for simplicity the equation
$$
u_{tt}-\2 u+u^5+u|u|^{p-2}=0\quad \roman{in}\e \Bbb
R^3\times[0,\infty[\tag5.1$$ 
with initial data (1.2). (5.1) may be approximated by 
equations
$$u_{tt}-\2
u+u^5\Bigl(1+
\min\bigl\{|u|^{p-6},k^{p-6}\bigr\}\Bigr)=0.\tag5.2$$
By the preceding, (5.2) admits global $C^2$-solutions 
$u^{(k)}$; moreover, as in
Rauch's Theorem 3.3 we may decompose
$$u^{(k)}(z_0)=u^{(0)}(z_0)+v^{(k)}(z_0),$$
where
$u^{(0)}$
solves the homogenous wave equation with initial data
$u_0,u_1$. Now, for suitable initial data, we obtain a 
uniform bound
$\bigl|u^{(0)}(z)\bigr|\le m_0\quad\hbox{for all}\e z\in 
R^3\times [0,\infty[;$
for instance, if
$u_0,u_1$
have compact support. By (2.3), moreover, if
$$
\bigl|u^{(k)}(z_0)\bigr|=\sup\limits
_{z\in K_0(z_0)}\bigl|u^{(k)}(z)\bigr|=m_k,
$$
we may estimate
$$
\split
m_k&=\bigl|u^{(k)}(z_0)\bigr|
\le\bigl|u^{(0)}(z_0)\bigr|+\bigl|v^{(k)}(z_0)\bigr|\\
&\le m_0+{m_k\over
4\pi}\il_{M_0(z_0)}{|u^{(k)}|^4\Bigl(1+
\min\bigl\{|u^{(k)}|^{p-6},k^{p-6}\bigr\}\Bigr)\over
t_0-t}\,do\\
&\le m_0+m_k{1+k^{p-6}\over 
4\pi}\il_{M_0(z_0)}{|u^{(k)}|^4\over t_0-t}\,do\\
&\le m_0+C\th m_k
k^{p-6}\Bigl(E\bigl(u(0)\bigr)+
E^{2/3}\bigl(u(0)\bigr)\Bigr)<2m_0,\endsplit
$$
if
$k=2m_0$
and if
$E\bigl(u(0)\bigr)$
is sufficiently small, depending on
$m_0$,
that is, on
$u_0$
and
$u_1$. Thus
$u=u^{(k)}$
solves (5.1). 

In particular, we obtain the following perturbation result:
%
\thm{Theorem 5.1}
Suppose
$u_0\in C^3(\bbbr^3)$,
$u_1\in C^2(\bbbr^3)$
have finite energy
$$
\il_{\bbbr^3}\left({|u_1|^2+|\nabla u_0|^2\over 2}+
{|u_0|^6\over 6}\right)dx <
\infty$$ and suppose the solution
$u^{(0)}$
to the homogeneous wave equation with initial data
$u_0,u_1$ is uniformly bounded. Then
there exists $\EPS_0>0$
such that for all $|\EPS|<\EPS_0$
the initial value problem for {\rm(5.1)} with data
$\EPS u_0$,
$\EPS u_1$
admits a global
$C^2$-solution.\ethm

However, ``in the large" the super-critical case appears 
to be completely open.


\Refs
\ref\no 1
\by P. Brenner
\paper On $L_p-L_{p'}$ estimates for the wave equation
\jour Math. Z. \vol 145 \yr 1975 \pages 251--254
\endref
\ref\no 2
\by P. Brenner and W. von Wahl
\paper Global classical solutions of non-linear wave 
equations
\jour Math. Z. \vol 176 \yr 1981 \pages 87--121\endref
\ref\no 3
\by F. E. Browder
\paper On nonlinear wave equations \jour Math. Z. \vol 80 
\yr 1962
\pages 249--264\endref
\ref\no 4
\by J. Ginibre and G. Velo
\paper The global Cauchy problem for the non-linear 
Klein-Gordon equation
\jour Math. Z. \vol 189 \yr 1985 \pages 487--505\endref
\ref\no 5
\bysame \paper Scattering theory in the energy space for 
a large class of
non-linear wave equations \jour Comm. Math. Phys. \vol 123
\yr 1989 \pages 535--573\endref
\ref\no 6
\by M. G. Grillakis \paper Regularity and asymptotic 
behaviour of the wave
equation with a critical nonlinearity \jour Ann. of Math. 
(2) (to appear)
\endref
\ref\no 7
\by F. John \paper Blow-up solutions to nonlinear wave 
equations in three
space dimensions \jour Manuscripta Math. \vol 28 \yr 1979
\pages 235--268\endref
\ref\no 8
\by K. J\"orgens \paper Das Anfangswertproblem im Grossen 
f\"ur
eine Klasse nicht-linearer Wellengleichungen \jour Math. 
Z. \vol
77 \yr 1961 \pages 295--308\endref
\ref\no 9
\by J. M. Lee and T. H. Parker \paper The Yamabe problem 
\jour Bull. Amer. Math.
Soc. (N.S.) \vol 17 \yr 1987 \pages 37--92\endref
\ref\no 10
\by J.-L. Lions \book Quelques m\'ethodes de r\'esolution 
des probl\`emes
aux limites non lin\'eaires \bookinfo Dunod, 
Gauthier-Villars, Paris
\yr 1969\endref
\ref\no 11 \by H. Pecher \paper Ein nichtlinearer 
Interpolationssatz und
seine Anwendung auf nichtlineare Wellengleichungen \jour 
Math. Z.
\vol 161 \yr 1978 \pages 9--40\endref
\ref\no 12
\by J. Rauch \book The $u^5$-Klein-Gordon equation {\rm 
(Brezis and
Lions, eds.)} \bookinfo Pitman Research Notes in Math. 
\publ no. 53 \pages
335--364\endref
\ref\no 13
\by L. I. Schiff \paper Nonlinear meson theory of nuclear 
forces
{\rm I} \jour Phys. Rev. \vol 84 \yr 1951 \pages 
1--9\endref
\ref\no 14
\by I. E. Segal \paper The global Cauchy problem for a 
relativistic
scalar field with power interaction \jour Bull. Soc. 
Math. France
\vol 91 \yr 1963 \pages 129--135\endref
\ref\no 15
\by J. Shatah \paper Weak solutions and the development 
of singularities
in the $SU(2)$ $\sigma$-model \jour Comm. Pure Appl. 
Math. \vol 41
\yr 1988 \pages 459--469\endref
\ref\no 16
\by W. Strauss \book Nonlinear wave equations \bookinfo 
CBMS
Lecture Notes, no. 73, Amer. Math. Soc., Providence, RI 
\yr 1989\endref
\ref\no 17
\by R. Strichartz \paper Restrictions of Fourier 
transforms to quadratic
surfaces and decay of solutions of wave equations \jour 
Duke Math. J.
\vol 44 \yr 1977 \pages 705--714\endref
\ref\no 18
\by M. Struwe \paper Globally regular solutions to the 
$u^5$-Klein-Gordon
equation \jour Ann. Sc. Norm. Sup. Pisa (Ser. 4) \vol 15
\yr 1988 \pages 495--513\endref
\ref\no 19
\by Y. Zheng \paper Concentration in sequences of 
solutions to the nonlinear
Klein-Gordon equation \jour preprint,  1989\endref
\ref\no 20
\by L. V. Kapitanski\u\i \paper The Cauchy problem for a 
semi-linear wave
equation, {\rm1989}\endref
\ref\no 21
\by J. Shatah and A. Tahvildar-Zadeh \paper Regularity of 
harmonic maps
from Minkowski space into rotationally symmetric 
manifolds \jour
Courant Institute, preprint,  1990\endref

\endRefs
\enddocument